\documentclass[11pt]{article}

\usepackage{empheq}
\usepackage{amsmath,amsfonts,amssymb,amsbsy,amscd}
\usepackage{mathrsfs,array}
\usepackage{subcaption}
\usepackage[a4paper, total={6in, 8in}]{geometry}
\usepackage{amsthm}

\newtheorem{theorem}{Theorem}[section]
\newtheorem{lemma}[theorem]{Lemma}
\newtheorem{remark}[theorem]{Remark}%

\usepackage{tikz}
\usetikzlibrary{arrows.meta}
\usetikzlibrary{3d}
\usetikzlibrary{shadings}
\usetikzlibrary{math}
\usepackage{adjustbox}\usepackage{multirow}
\usepackage{upgreek}
\usetikzlibrary{intersections, arrows, automata, backgrounds, calendar, chains, matrix, mindmap, patterns, petri, shadows, shapes.geometric, shapes.misc, spy, trees, shapes}

\definecolor{darkgreen}{rgb}{0, 0.7, 0}
\definecolor{orange}{rgb}{0.98, 0.6, 0.01}
 	\definecolor{napiergreen}{rgb}{0.16, 0.5, 0.0}

\usepackage{upgreek}

\newcolumntype{R}[2]{%
    >{\adjustbox{angle=#1,lap=\width-(#2)}\bgroup}%
    l%
    <{\egroup}%
}
\newcommand*\rot{\multicolumn{1}{R{45}{1em}}}
\makeatletter
\newcommand{\thickhline}{%
    \noalign {\ifnum 0=`}\fi \hrule height 1pt
    \futurelet \reserved@a \@xhline
}
\newcolumntype{"}{@{\hskip\tabcolsep\vrule width 1pt\hskip\tabcolsep}}
\makeatother

\usepackage{pifont}
\newcommand{\cmark}{\ding{51}}%
\newcommand{\xmark}{\ding{55}}%
\newcommand{\smark}{\ding{72}}%

\newcommand\xqed[1]{%
  \leavevmode\unskip\penalty9999 \hbox{}\nobreak\hfill
  \quad\hbox{#1}}
\newcommand\demo{\xqed{$\triangle$}}

\usepackage{hyperref}
\hypersetup{
    colorlinks=true,
    linkcolor=blue,
    urlcolor=red,
    linktoc=all,
    allcolors=blue
}

\usepackage[nameinlink,capitalise,noabbrev]{cleveref}
\makeatletter
\newcommand{\myref}[1]{\cref{#1}\mynameref{#1}{\csname r@#1\endcsname}}
\newcommand{\Myref}[1]{\Cref{#1}\mynameref{#1}{\csname r@#1\endcsname}}

\usepackage{accents}
\newlength{\dhatheight}
\newcommand{\doublehat}[1]{%
    \settoheight{\dhatheight}{\ensuremath{\hat{#1}}}%
    \addtolength{\dhatheight}{-0.20ex}%
    \hat{\vphantom{\rule{1pt}{\dhatheight}}%
    \smash{\hat{#1}}}}
\newcommand{\doublecheck}[1]{%
    \settoheight{\dhatheight}{\ensuremath{\check{#1}}}%
    \addtolength{\dhatheight}{-0.10ex}%
    \check{\vphantom{\rule{1pt}{\dhatheight}}%
    \smash{\check{#1}}}}

\allowdisplaybreaks

\title{A unified framework for Navier-Stokes Cahn-Hilliard models with non-matching densities}

\author{M.F.P. ten Eikelder$^{\dag,\ddag,}$\thanks{Corresponding author. e-mail: \texttt{marco.eikelder@tu-darmstadt.de}}
\and K.G. van der Zee$^\S$
\and I. Akkerman$^\P$
\and D. Schillinger$^\dag$
}
\date{%
    $^\dag$Institute for Mechanics, Computational Mechanics Group, Technical University of Darmstadt\\
    $^\ddag$Institute of Mechanics and Computational Mechanics, Leibniz University Hannover\\
    $^\S$School of Mathematical Sciences, University of Nottingham\\
    $^\P$Department of Mechanical, Maritime and Materials Engineering, Delft University of Technology    
}

\def\B#1{\mbox{\boldmath{$#1$}}}

\newcommand{\divg}{{\rm div}}

\newcommand{\nn}{\nonumber}

\newcommand{\bu}{\mathbf{u}}
\newcommand{\bD}{\mathbf{D}}

\newcommand{\jv}{[\![\bv]\!]}
\newcommand{\jrho}{[\![\rho]\!]}
\newcommand{\arho}{\left\{\rho\right\}}
\newcommand{\bw}{\B{w}}

\newcommand{\bv}{\mathbf{v}}

\newcommand{\bx}{\mathbf{x}}

\newcommand{\bg}{\B{g}}
\newcommand{\bJ}{\mathbf{J}}
\newcommand{\bh}{\mathbf{h}}
\newcommand{\bomega}{\boldsymbol{\omega}}
\newcommand{\bchi}{\boldsymbol{\chi}}
\newcommand{\bnu}{\boldsymbol{\nu}}
\newcommand{\trho}{\tilde{\rho}}
\newcommand{\bpi}{\boldsymbol{\pi}}

\def\be{\begin{equation}}
\def\ee{\end{equation}}
\def\ba{\begin{array}}
\def\ea{\end{array}}
\def\bea{\begin{eqnarray}}
\def\eea{\end{eqnarray}}
\def\beas{\begin{eqnarray*}}
\def\eeas{\end{eqnarray*}}
\newcommand{\bseq}{\begin{subequations}}
\newcommand{\eseq}{\end{subequations}}

\begin{document}

\maketitle

\begin{abstract}
Over the last decades, many diffuse-interface Navier-Stokes Cahn-Hilliard models with non-matching densities have appeared in the literature. These models claim to describe the same physical phenomena, yet they are distinct from one another. The overarching objective of this work is to bring all of these models together by laying down a unified framework of Navier-Stokes Cahn-Hilliard models with non-zero mass fluxes. Our development is based on three unifying principles: (1) there is only one system of balance laws based on continuum mixture theory that describes the physical model, (2) there is only one natural energy-dissipation law that leads to quasi-incompressible Navier-Stokes Cahn-Hilliard models, (3) variations between the models only appear in the constitutive choices. The framework presented in this work now completes the fundamental exploration of alternate non-matching density Navier-Stokes Cahn-Hilliard models that utilize a single momentum equation for the mixture velocity, but leaves open room for further sophistication in the energy functional and constitutive dependence.
\end{abstract}

\noindent{\small{\textbf{Key words}. Navier-Stokes Cahn-Hilliard equations, phase-field models, incompressible two-phase flow, mixture theory, thermodynamic consistency.}}\\


\noindent{\small{\textbf{AMS Subject Classification}:  Primary: 76T99, Secondary: 35Q30, 35Q35, 35R35, 76D05, 76D45, 80A99}}

\section{Introduction}\label{sec: intro}
Phase-field models have emerged as a powerful tool for describing interface problems that appear in various fields in science. These models are also termed diffuse-interface models due to the smooth representation of the interface  \cite{anderson1998diffuse}. The theory of general diffuse-interface models in solid and fluid mechanics has been presented in   \cite{oden2010general}. Phase-field models are typically equipped with a thermodynamical framework  \cite{malek2008thermodynamic,muller2001thermodynamics,penrose1990thermodynamically,souvcek2014natural} and stability properties  \cite{elliott1990global}. A key factor in the success of phase-field models is their ability to be directly applied in computer simulations  \cite{gomez2018computational} ranging from phase transitions in fluids \cite{liu2015liquid} to fracture mechanics \cite{borden2012phase}.

In the context of free-surface fluid mechanics, other popular methodology are the volume-of-fluid methods \cite{Hirt_Nichols_81} and level-set methods \cite{osher2006level,sethian1999level}. Level-set methods are popular for incompressible flows \cite{akkerman2011isogeometric,carrica2007unsteady,ten2021novel} and volume-of-fluid methods are employed for both incompressible \cite{gueyffier1999volume,park2009volume} and compressible flows \cite{ten2017acoustic,murrone2005five,saurel1999multiphase}. The main distinguishing feature of a phase-field model, as compared with volume-of-fluid and level-set models, is that the interface is determined by a physical model with a fixed interface width. In computations with the volume-of-fluid approach often an interpolation technique or a compression algorithm is applied near the interface, whereas level-set methods typically use a redistancing algorithm to maintain a stable interface. Apart from the thermodynamical structure of the underlying phase-field formulations, redundancy of such algorithms is often the most important practical advantage of phase-field simulations.

The phase-field model that describes incompressible isothermal two constituent flows with non-matching densities is the Navier-Stokes Cahn-Hilliard (NSCH) model. Over the years many NSCH models with non-matching densities have been proposed with distinct fundamental variables for the velocity and phase-field. It is the purpose of the current paper to present a unified framework for non-matching density NSCH models which is invariant to the choice of fundamental variables. The framework encompasses variations of existing NSCH models.

\subsection{Historical overview}
The first coupling between the Navier-Stokes equations, describing viscous fluid flow, and the Cahn-Hilliard equation, describing spinoidal decomposition, has been established by Hohenberg and Halperin \cite{hohenberg1977theory}. They proposed the system, referred to as \textit{model H}, that reads:
\begin{subequations}\label{eq: model H}
  \begin{align}
     \rho \partial_t \bu + \rho(\bu \cdot \nabla \bu)\bu - \divg \left(2 \nu(c) \bD \right) + \nabla p =&~ -\sigma\epsilon \divg (\nabla c \otimes \nabla c),\label{eq: model H: mom eq}\\
     \divg \bu =&~ 0,\label{eq: model H: cont eq}\\
     \partial_t c + \bu \cdot \nabla c =&~ \divg (m(c) \nabla \mu),\label{eq: model H: CH eq}\\
     \mu =&~ \sigma\epsilon^{-1} \varphi'(c) - \sigma\epsilon \Delta c, \label{eq: model H: chemical pot}
  \end{align}
\end{subequations}
in domain $\Omega \subset \mathbb{R}^d$, with dimension $d = 2,3$, both open and bounded that is occupied with two constituents $j=1,2$. We adopt the standard notation, where $\mathbf{x}$ represents the (Eulerian) position vector, $t$ is the time, $\partial_t$ the partial time-derivative, $\nabla$ the gradient, ${\rm div}$ the divergence (defined as $({\rm div} \mathbf{A})_i:= \partial A_{ik}/\partial x_k$) and $\Delta$ the Laplace operator. Here $\bu$ is the mean velocity, $\rho$ is the density, $p$ is the pressure and $c$ is the so-called order parameter representing a concentration related quantity. Moreover $\bD$ represents the symmetric gradient of the mean velocity, $\bD = \frac{1}{2}\left(\nabla \bu + (\nabla \bu)^T\right)$, $\nu(c)$ is the concentration dependent dynamic viscosity of the mixture, $\sigma$ is the surface tension coefficient and $\epsilon$ is an interface thickness parameter. The surface tension coefficient is assumed constant, i.e. Marangoni-type effects are precluded in this model. The quantity $\varphi=\varphi(c)$ is the homogeneous free energy and $m=m(c)\geq 0$ is the mobility. In model \eqref{eq: model H}, equation \eqref{eq: model H: mom eq} describes the linear momentum equation in which the right-hand side, $-\sigma\epsilon \divg (\nabla c \otimes \nabla c)$, represents a contribution of the stress tensor that models capillary forces due to surface tension. Equation \eqref{eq: model H: cont eq} dictates that the mean velocity is divergence-free. Lastly, the Cahn-Hilliard equation \eqref{eq: model H: CH eq} describes the evolution of the order parameter in which the right-hand side represents the divergence of a diffusive flux, and the variable $\mu$ given by \eqref{eq: model H: chemical pot} is often referred to as the chemical potential.

The major limitation of model H is its assumption of constant density, i.e. the density of the mixture as well as the density of the individual constituents is constant. As such, this precludes the applicability of the model to problems with large density ratios (e.g. water-air flows). Model H has been used to study the fluid behavior at critical points of single and binary fluids. The derivation of the model initially relied on phenomenological arguments and a rigorous derivation in the framework of rational continuum mechanics was absent until 1996.

In that year Gurtin \cite{gurtinmodel} provided this missing derivation. The core ideas in this derivation are (i) the separation of balance laws from constitutive relations, (ii) the usage of microforce balance laws, and (iii) the selection of constitutive equations compatible with an energy-dissipation law referred to as the second law of thermodynamics (the so-called Coleman-Noll procedure). Gurtin introduces the balance of mass per constituent while he directly presents a single momentum equation for the mixture. In the same year Gurtin applied this procedure to present the Ginzburg-Landau and Cahn-Hilliard equations with a derivation in a rational continuum mechanics framework \cite{gurtin1996generalized}. Several extensions of the microforce continuum framework of Gurtin have been proposed, including a phase-field gradient theory for enriched continua \cite{espath2021phase}.

One of the first efforts of extending model H to the case of non-matching densities is the work of Lowengrub and Truskinovsky \cite{lowengrub1998quasi}. Lowengrub and Truskinovsky present the following \textit{quasi-incompressible model}:
\begin{subequations}\label{eq: intro model lowengrub}
  \begin{empheq}[left=\empheqlbrace]{align}
   \partial_t (\rho \bv) + {\rm div} \left( \rho \bv\otimes \bv \right) + \nabla p 
    - {\rm div} \left(   \nu(c) (2\mathbf{D}+\lambda(c)({\rm div}\bv) \mathbf{I}) \right)&\nn\\
    + \sigma\epsilon {\rm div} \left( \rho \nabla c \otimes \nabla c \right) &=~ 0, \label{eq: intro model lowengrub: mom}\\
 \partial_t \rho + {\rm div} (\rho \bv) &=~ 0, \label{eq: intro model lowengrub: cont} \\
  \rho (\partial_t c + \bv \cdot \nabla c) -\divg (m \nabla \mu) &=~0,\label{eq: intro model lowengrub: PF}\\
  \mu -\frac{\sigma}{\epsilon} \dfrac{\partial \varphi}{\partial c} + \dfrac{\sigma\epsilon}{\rho}{\rm div}\left( \rho  \nabla c \right) +\dfrac{p}{\rho^2} \dfrac{\partial \rho}{\partial c}&=~0,\label{eq: intro model lowengrub: mu}
  \end{empheq}
\end{subequations}
with free energy density per unit volume $\rho \sigma\left( \epsilon^{-1}\varphi(c)+\epsilon|\nabla c|^2/2 \right)$ and $m\geq 0$ a non-degenerate mobility (the mobility is called \textit{non-degenerate} if it is constant and \textit{degenerate} if it vanishes in the single-fluid regime). The \textit{mixture density} $\rho$ is  $\rho = \trho_1 + \trho_2$ and the \textit{mass-averaged velocity} $\bv$ is given by $\rho \bv = \trho_1 \bv_1 + \trho_2 \bv_2$ where $\trho_{j}=\rho c_{j}$ represents the variable density of constituent $j$ with $c_{j}$ the concentration of constituent $j$. In comparison with model H, this model has two key differences. First, the incompressibility constraint \eqref{eq: model H: cont eq} is replaced by the conservation of mass of the mixture \eqref{eq: intro model lowengrub: cont}. Expanding the divergence in \eqref{eq: intro model lowengrub: cont} reveals the quasi-incompressible nature of the model:
\begin{align}
     {\rm div} \bv + \frac{1}{\rho} \dfrac{\partial \rho}{\partial c} (\partial_t c + \bv\cdot \nabla c) &=~ 0.
\end{align}
Second, the evolution equation of the order parameter \eqref{eq: intro model lowengrub: PF} explicitly contains a density which distinguishes it from the evolution equation \eqref{eq: model H: CH eq}. It seems challenging to design algorithms that solve \eqref{eq: intro model lowengrub} and often a simplified version of the model is employed, see e.g.  \cite{lee2002modeling}. The challenge is often linked to the quasi-incompressible nature of the model \cite{abels2012thermodynamically}. We note that recently a numerical method for a reformulation of the model \eqref{eq: intro model lowengrub} has been proposed in  \cite{guo2014numerical}. 

Noteworthy alternative models have been proposed by Boyer \cite{boyer2002theoretical} and Ding et al. \cite{ding2007diffuse}. Both models deviate from the model of Lowengrub and Truskinovsky \cite{lowengrub1998quasi} in that they use a mean velocity that is fully incompressible:
\begin{align}
    \divg \bu = 0.
\end{align}
This mean velocity, defined as $\bu = \phi_1 \bv_1 + \phi_2 \bv_2$ where $\phi_{j}$ is the volume fraction of constituent $j$, is referred to as the \textit{volume-averaged velocity}. Actually the model Ding et al. \cite{ding2007diffuse} is, apart from the variable density, identical to the model H (model \eqref{eq: model H}). Analogously to Gurtin \cite{gurtinmodel} and Lowengrub and Truskinovsky \cite{lowengrub1998quasi}, the starting point of Ding et al. \cite{ding2007diffuse} consists of the individual mass balance equations of the constituents and the momentum equation of the mixture. In contrast, the point of departure adopted in the model of Boyer \cite{boyer2002theoretical} is composed of individual mass and momentum balance equations, which leads to a different form of the momentum and phase equations.

Unfortunately the models of Boyer \cite{boyer2002theoretical} and Ding et al. \cite{ding2007diffuse} are not presented with some (approximate) form of the second law of thermodynamics. This observation has led to the development of the model of Abels et al. \cite{abels2012thermodynamically} which reads:
\begin{subequations}\label{eq: intro model abels}
  \begin{empheq}[left=\empheqlbrace]{align}
   \partial_t (\rho \bu) + {\rm div} \left( \rho \bu\otimes \bu \right) + {\rm div}\left( \bu \otimes \frac{\rho_1-\rho_2}{2}m(\phi)\nabla \mu \right) + \nabla p 
    &\nn\\- {\rm div} \left(  2 \nu(\phi) \mathbf{D} \right) +\sigma\epsilon {\rm div} \left( \nabla \phi \otimes \nabla \phi \right)&=~0, \label{eq: intro model abels: mom}\\
{\rm div} \bu &=~ 0, \label{eq: intro model abels: cont} \\
 \partial_t \phi + \bu \cdot \nabla \phi -\divg (m(\phi) \nabla \mu) &=~0,\label{eq: intro model abels: PF}\\
  \mu -\frac{\sigma}{\epsilon} \dfrac{\partial \varphi}{\partial \phi}+\sigma\epsilon\Delta \phi&=~ 0 ,\label{eq: intro model abels: mu}
  \end{empheq}
\end{subequations}
and is presented with some form of the second law of thermodynamics. This form of the second law is  distinct from what appears in Lowengrub and Truskinovsky \cite{lowengrub1998quasi} in the sense that a different kinetic energy is used. This model also uses the volume-averaged mean velocity for the velocity of the mixture and makes use of the difference of volume fractions as order parameter. The distinguishing feature of \eqref{eq: intro model abels} is the, somehow surprising, additional convective term in the momentum equation.

Since then several quasi-incompressible models that employ the difference of volume fractions as order parameter have been proposed. Noteworthy contributions are the work of Shen et al. \cite{shen2013mass},  Aki et al. \cite{aki2014quasi} and Shokrpour Roudbari et al. \cite{shokrpour2018diffuse}, which are all presented with an approximate form of the second law. The models \cite{aki2014quasi} and  \cite{shokrpour2018diffuse} have been derived in a similar way, namely using balance laws of the individual constituents and the Coleman-Noll procedure. The derivation of the model of Shen et al. follows different considerations. An important contribution in the work of Shokrpour Roudbari et al. \cite{shokrpour2018diffuse} is the observation that, up to the definition of the mobility and the definition of mass fluxes, the models of Aki et al. \cite{aki2014quasi},  Shen et al. \cite{shen2013mass} and Shokrpour Roudbari et al. \cite{shokrpour2018diffuse} are all equivalent.

\subsection{Objective and main results}
It is obvious that the literature on NSCH models is divided. The various proposed models (almost) all aim to represent the same physical behavior but are either clearly different or are presented in a (sometimes seemingly) different form. The four key flavors are:
\begin{enumerate}
    \item the mixture velocity: mass-averaged velocity ($\bv$) or a volume-averaged velocity ($\bu$)
    \item the order parameter: based on volume fractions ($\phi$) or concentration ($c$)
    \item the type of the free energy: volume-measure-based ($\Psi$) or a mass-measure-based ($\psi$)
    \item the mobility: non-degenerate (i.e. constant) or degenerate (i.e. it vanishes in the single-fluid regime) 
\end{enumerate}
We provide a summary of several existing models in \cref{table: overview models}.
\begin{center}
\begin{table}[h!]
{\small
\begin{tabular}{m{14em}m{1.8em}m{1.8em}m{1.8em}m{8.0em}m{1.8em}}
\textbf{Model}                      & \rot{\textbf{Velocity}} & \rot{\textbf{Order parameter}} & \rot{\textbf{Free energy}} & \rot{\textbf{Mobility}}     & \rot{\textbf{Energy law}}   \\[6pt] \thickhline\\[-4pt]
Abels et al. \cite{abels2012thermodynamically}  & $\bu$         & $\phi$          & $\Psi$      & non-deg./degen.      & {\color{orange}\smark}                 \\[6pt] 
Aki et al. \cite{aki2014quasi}                  & $\bv$         & $\phi$          & $\Psi$      & non-deg.    & {\color{darkgreen}\cmark}              \\[6pt] 
Boyer \cite{boyer2002theoretical}                    & $\bu$         & $\phi$          & $\Psi$      & degen.   & {\color{red}\xmark}               \\[6pt] 
Ding et al. \cite{ding2007diffuse}                & $\bu$         & $\phi$          & $\Psi$      & degen.   & {\color{red}\xmark}                   \\[6pt] 
Lowengrub and Truskinovsky \cite{lowengrub1998quasi} & $\bv$         & $c$             & $\psi$      & non-deg. & {\color{darkgreen}\cmark}                   \\[6pt] 
Shen et al. \cite{shen2013mass}               & $\bv$         & $\phi$          & $\Psi$      & non-deg. & {\color{darkgreen}\cmark}                   \\[6pt] 
Shokrpour Roudbari et al. \cite{shokrpour2018diffuse}  & $\bv$         & $\phi$          & $\Psi$      & non-deg. & {\color{darkgreen}\cmark}              \\[6pt] \hline
\end{tabular}}
\caption{Overview of various NSCH models. The columns indicate the mean velocity (either mass-averaged $\bv$ or volume-averaged $\bu$), the order parameter (either volume fraction based $\phi$ or concentration based $c$), the type of the free energy (either volume-measure-based $\Psi$ or mass-measure-based $\psi$), the type of the mobility (either degenerate or non-degenerate) and lastly whether the model is equipped with an energy-dissipation law. The model of Abels et al. is presented with both non-degenerate and degenerate mobilities. Additionally, the energy-dissipation law of that model differs from that of the other models, due to the usage of the volume-averaged velocity, and we indicate this by the symbol {\color{orange}\smark}. }
\label{table: overview models}
\end{table}
\end{center}
Even though each of these works has provided new and useful insights into diffuse-interface modeling and some provide elegant and physically sound derivations, the current status is far from optimal. We have two main objections:
\begin{enumerate}
    \item the systems of balance laws of the various models are distinct before constitutive choices have been applied,
    \item the energy-dissipation laws of the models are not identical.
\end{enumerate}
Since all the models represent the same physics, we lay down the three unifying principles:
\begin{enumerate}
    \item there is only one system of continuum mechanics balance laws that describes the physical model,
    \item there is only one natural energy-dissipation law that leads to quasi-incompressible NSCH models,
    \item variations between the models only appear in the constitutive choices.
\end{enumerate}
The main objective of this work is to lay down a \textit{unified framework of incompressible NSCH models with non-zero mass fluxes} on the basis of these three unifying principles. In particular we establish one incompressible NSCH system of balance laws and show that many alternate forms are connected via variable transformations. Two (equivalent) formulations that result after the constitutive choices are (i) a formulation in terms of the mass-averaged velocity $\bv$:
\begin{subequations}\label{eq: model v with const mod intro}
  \begin{empheq}[left=\empheqlbrace]{align}
   \partial_t (\rho \bv) + {\rm div} \left( \rho \bv\otimes \bv \right) + \nabla p + {\rm div} \left( \nabla \phi \otimes \dfrac{\partial  \Psi}{\partial \nabla \phi} + (\mu\phi-\Psi)\mathbf{I} \right) & \nn\\
    - {\rm div} \left(   \nu (2\mathbf{D}+\lambda({\rm div}\bv) \mathbf{I}) \right)-\rho\mathbf{g} &=~ 0, \label{eq: model v with const mod: mom intro}\\
 \partial_t \rho + {\rm div}(\rho \bv) &=~ 0, \label{eq: model v with const mod: cont intro} \\
   \partial_t \phi + {\rm div}(\phi \bv) - {\rm div} \left(\mathbf{M}^v\nabla (\mu+\alpha p)\right) +\zeta m (\mu + \alpha p) &=~0,\label{eq: model v with const mod: PF intro}\\
  \mu - \dfrac{\partial \Psi}{\partial \phi}+{\rm div} \left(  \dfrac{\partial \Psi}{\partial \nabla \phi} \right)&=~0,
  \end{empheq}
\end{subequations}
and (ii) a formulation in terms of the volume-averaged velocity $\bu$:
\begin{subequations}\label{eq: model u with const mod intro}
  \begin{empheq}[left=\empheqlbrace]{align}
    \partial_t (\rho \bu + \tilde{\bJ}^u) + {\rm div} \left( \rho \bu\otimes \bu + \tilde{\bJ}^u\otimes \bu + \bu\otimes \tilde{\bJ}^u + \frac{1}{\rho}\tilde{\bJ}^u\otimes\tilde{\bJ}^u \right)&\nn\\
    + \nabla p + {\rm div} \left( \nabla \phi \otimes \dfrac{\partial \Psi}{\partial \nabla \phi} + (\mu\phi-\Psi)\mathbf{I} \right)\nn\\
    - {\rm div} \left(   \nu \left(2\nabla^s\left(\bu + \rho^{-1}\tilde{\bJ}^u\right)+\lambda\left({\rm div}\left(\bu +\rho^{-1}\tilde{\bJ}^u\right)\right) \mathbf{I}\right) \right)-\rho\mathbf{g}&=~ 0, \label{eq: model u with const mod: mom intro}\\
 {\rm div} \bu -\beta \gamma&=~ 0, \label{eq: model u with const mod: cont} \\
  \partial_t \phi + \bu\cdot \nabla\phi- {\rm div} \left(\mathbf{M}^u\nabla (\mu+\alpha p)\right)+\dfrac{\rho}{2\rho_1\rho_2} m (\mu + \alpha p)&~=0\label{eq: model u with const mod: PF intro},\\
 \mu - \dfrac{\partial\Psi}{\partial \phi}+{\rm div} \left(  \dfrac{\partial \Psi}{\partial \nabla \phi} \right)&=~0,
  \end{empheq}
\end{subequations}
with $\tilde{\bJ}^u=-(\rho_1-\rho_2)\mathbf{M}^u\nabla(\mu+\alpha p)/2$. Here $\mathbf{M}^v=\mathbf{M}^v(\phi,\nabla \phi, \mu, \nabla \mu, p)$, $\mathbf{M}^u(\phi,\nabla \phi, \mu, \nabla \mu, p)$ and $m=m(\phi, \mu, p)$ are degenerate mobilities, $\nu=\nu(\phi)$ is the dynamic viscosity, $p$ is the pressure, $\bg$ is the gravitational acceleration, $\rho_1$ and $\rho_2$ are the constant specific densities of the constituents and we have introduced the constants $\alpha = (\rho_2-\rho_1)/(\rho_1+\rho_2)$, $\beta=(\rho_2-\rho_1)/(2\rho_1\rho_2)$ and  $\zeta=(\rho_1+\rho_2)/(2\rho_1\rho_2)$. We refer the reader for details on the specific choices and scaling of the mobility parameters, and free energy to Abels et al. \cite{abels2012thermodynamically},  Aki et al. \cite{aki2014quasi} and  Gurtin \cite{gurtinmodel}. Moreover we show that (i) many existing models are identical (up to the definition of the mobility), (ii) existing volume-averaged velocity based models are inconsistent with the mixture theory framework, however have a consistent rectification. As a side observation it turns out that, within our framework, models with a non-degenerate mobility are incompatible in the single-fluid region (see \cref{rmk: incomp mob}).

\subsection{Plan of the paper}
The structure of the remainder of the paper is as follows. In \cref{sec: mix theory} the system of balance laws is established. We present the framework of balance laws for binary mixtures of incompressible viscous fluids. Starting from constituent balance laws, we use mixture theory to derive the balance laws of the mixtures. Additionally, we present some new and important identities and evolution equations. In \cref{sec: 2nd law} we perform constitutive modeling via the Coleman-Noll procedure. Here we highlight the modeling assumptions of the NSCH model from the viewpoint of mixture theory. Additionally, we show that applying the Coleman-Noll procedure to alternative derivations provides the same modeling restriction. In \cref{sec: Relation to existing Navier-Stokes Cahn-Hilliard models} we discuss the relation of the novel model to existing NSCH models. Finally, in \cref{sec: discussion} we summarize our findings and present some possible further research directions.

\section{Mixture theory}\label{sec: mix theory}
In this section we lay down the mixture theory as well as the necessary definitions. In this work we focus on incompressible isothermal constituents, which we will specify in \cref{sec: prelim}. We restrict ourselves to the case of binary mixtures for the sake of simplicity and note that the multi-component case is a straightforward extension. This section is fully compatible (i.e. no approximations are introduced) with the mixture theory metaphysical principles proposed by Truesdell \cite{truesdell1984historical}:
\begin{enumerate}
    \item All properties of the mixture must be
mathematical consequences of properties of the
constituents.
\item So as to describe the motion of a constituent, we may in imagination isolate it from the rest of the mixture, provided we allow properly for the actions of the other constituents upon it.
\item The motion of the mixture is governed
by the same equations as in a single body.
\end{enumerate}
\cref{sec: prelim} provides the necessary (kinematic) definitions. Then in \cref{sec: BL single constituents,sec: BL mixtures} we introduce the balance laws of the individual constituents and subsequently, relying on metaphysical principles of  \cite{truesdell1960classical}, balance laws of the mixtures. Finally in \cref{sec: alternative formulations} we present alternative, but equivalent, formulations of the balance laws, diffusive fluxes and stress tensors.

\subsection{Preliminaries}\label{sec: prelim}
In domain $\Omega$ we denote with $\mathbf{X}_{j}$ the position of a particle of constituent $j$ in the Lagrangian (reference) configuration. Denoting the position of the mixture in the Eulerian frame by $\mathbf{x}$, we identify its relation to the initial configurations of the constituents with an invertible deformation map $\bchi_{j}$:
\begin{align}
    \mathbf{x} := \bchi_{j}(\mathbf{X}_{j},t). 
\end{align}
The velocity of constituent $j$ is given by
\begin{align}
    \bv_{j}(\bx,t):=\partial_t\bchi_{j}(\mathbf{X}_{j},t)=\partial_t\bchi_{j}(\bchi_{j}^{-1}(\bx,t),t).
\end{align}
Consider now an arbitrary control volume $V \subset \Omega$ around spatial position $\mathbf{x}$ in the mixture that contains the masses $M_{j}=M_{j}(V)$ of the constituents $j=1,2$ at time $t$. We denote the total mass in $V$ as $M=M(V)=\sum_{j}M_{j}(V)$. We define the partial mass density $\tilde{\rho}_{j}$ as the mass of constituent $j$ per unit volume of the mixture as:
\begin{align}\label{eq: def trhoA}
  \tilde{\rho}_{j}(\bx,t) := \displaystyle\lim_{|V| \rightarrow 0} \dfrac{M_{j}(V)}{|V|},
\end{align}
where $|V|$ denotes the measure of control volume $V$. The mass density of the total mixture at position $\bx$ and time $t$ is now defined as the sum of the partial mass densities of the constituents:
\begin{align}\label{eq: def rho}
\rho(\bx,t):=\displaystyle\sum_{j}\tilde{\rho}_{j}(\bx,t).
\end{align}
Next, we introduce $V_{j} \subset V$ as the control volume occupied by constituent $j$. We define the volume fraction of constituent $j$ as:
\begin{align}\label{eq: def phi}
  \phi_{j}(\bx,t) := \displaystyle\lim_{|V| \rightarrow 0} \dfrac{|V_{j}|}{|V|},
\end{align}
where $|V_{j}|$ denotes the measure of $V_{j}$.
We assume that
\begin{align}\label{eq: sum phi}
    \displaystyle\sum_{j} \phi_{j}= 1,
\end{align}
and thus exclude the existence of interstitial voids.
Additionally we introduce the concentration of constituent $j$ as:
\begin{align}\label{eq: def conc}
  c_{j}(\bx,t) := \displaystyle\lim_{|V| \rightarrow 0} \dfrac{M_{j}(V)}{M(V)}.
\end{align}
The above definitions \eqref{eq: def trhoA}, \eqref{eq: def rho} and \eqref{eq: def conc} imply the relation:
\begin{align}
  \tilde{\rho}_{j}(\bx,t)=\rho(\bx,t)c_{j}(\bx,t).
\end{align}
Next, we define the specific mass density $\rho_{j}$ as the mass of constituent $j$ per volume occupied by that constituent as:
\begin{align}\label{eq: def rhoA}
  \rho_{j}(\bx,t) := \displaystyle\lim_{|V_{j}| \rightarrow 0} \dfrac{M_{j}(V)}{|V_{j}|}.
\end{align}
We assume $\rho_{j}(\bx,t)>0$. From the definitions \eqref{eq: def trhoA}, \eqref{eq: def phi} and \eqref{eq: def rhoA} we deduce:
\begin{align}
  \tilde{\rho}_{j}(\bx,t) = \rho_{j}(\bx,t)\phi_{j}(\bx,t).
\end{align}
Additionally we have the relation:
\begin{align}\label{eq: sum c}
    \displaystyle\sum_{j} c_{j}= 1.
\end{align}
In general the specific mass densities $\rho_j$ may vary due to compressibility and thermal effects. In this paper we restrict to constant specific mass densities $\rho_{j}$ representing incompressible isothermal constituents. The momentum associated with constituent $j$ is:
\begin{align}
    \mathbf{m}_{j}(\bx,t) := \trho_{j}(\bx,t) \bv_{j}(\bx,t).
\end{align}
The momentum of the mixture is the sum of that of the individual constituents:
\begin{align}\label{eq: momentum mixture}
    \mathbf{m}(\bx,t) := \displaystyle\sum_{j} \mathbf{m}_{j}(\bx,t).
\end{align}
The \textit{mixture velocity} $\bv$ is a \textit{mass-averaged velocity} or \textit{barycentric velocity} and is identified via the relation:
\begin{align}\label{eq: mix velo}
    \mathbf{m}(\bx,t) = \rho(\bx,t) \bv(\bx,t).
\end{align}
The \textit{diffusion velocity} or \textit{peculiar velocity} of constituent $j$ is the constituent velocity relative to the gross motion of the mixture:
\begin{align}\label{eq: def bwj}
    \bw_{j}(\bx,t):=\bv_{j}(\bx,t)-\bv(\bx,t).
\end{align}
An immediate consequence is the observation that the superposition of the momenta relative to the gross motion of the mixture vanishes: 
\begin{align}\label{eq: rel gross motion zero}
    \displaystyle\sum_{j} \trho_{j} \bw_{j} = \displaystyle\sum_{j} \trho_{j} \bv_{j} - \displaystyle\sum_{j} \trho_{j} \bv = 0.
\end{align}
From \eqref{eq: rel gross motion zero} we deduce the identity:
\begin{align}\label{eq: sum diffusion velocity 2}
  \displaystyle\sum_{j} c_{j} \bw_{j} = 0,
\end{align}
that will be employed later in this section.
We now introduce two different material derivatives, one that follows the individual motion of constituent $j$ and one that follows the mean motion, respectively given by:
\begin{subequations}
\begin{align}
    \grave{\uppsi}=&~ \partial_t \uppsi + \bv_{j}\cdot \nabla \uppsi,\\
    \dot{\uppsi} =&~ \partial_t \uppsi + \bv\cdot \nabla \uppsi.\label{eq: mat der}
\end{align}
\end{subequations}
We define the jump and average of a constituent-related (vector-valued) quantity as:
\begin{subequations}
\begin{align}
    [\![\bw]\!]:=&~\frac{1}{2}(\bw_1-\bw_2),\\
    \left\{\bw\right\}:=&~\frac{1}{2}(\bw_1+\bw_2),
\end{align}
\end{subequations}
respectively, where the subscripts refer to the constituent numbers. Lastly, we introduce the constants:
\begin{subequations}\label{eq: def alpha beta}
\begin{align}
    \alpha :=&~ \frac{\rho_2-\rho_1}{\rho_1+\rho_2} =- \frac{\jrho}{\arho}, \label{eq: def alpha}\\
    \beta :=&~ \dfrac{\rho_2-\rho_1}{2\rho_1\rho_2}= -\dfrac{\jrho}{\rho_1\rho_2},\label{eq: definition beta}\\
    \zeta :=&~ \dfrac{\rho_1+\rho_2}{2\rho_1\rho_2}= \dfrac{\arho}{\rho_1\rho_2}.\label{eq: definition zeta}
\end{align}
\end{subequations}

\subsection{Balance laws of single constituents}\label{sec: BL single constituents}
Denoting by $\gamma_{j}$ the mass supply of constituent $j$ due to reaction, we introduce the local evolution equation of the mass of constituent $j$:
\begin{align}\label{eq: local mass balance constituent j}
    \partial_t \tilde{\rho}_{j} + {\rm div}(\tilde{\rho}_{j} \bv_{j})= \gamma_{j}.
\end{align}
The associated convective form is:
\begin{align}\label{eq: local mass balance constituent j conv}
   \grave{\trho}_j + \trho_{j}\divg \bv_j= \gamma_{j}.
\end{align}
We assume that mass fluxes $\gamma_{j}$ vanish in the single fluid region (i.e. when $\phi_{j}=\pm1$). As a consequence, the mass balance in the single fluid region reads:
\begin{align}
    \partial_t \rho_{j} + {\rm div}(\rho_{j} \bv_{j})= 0, \quad \text{for }\phi_{j}=1.
\end{align}
Since the specific constituent densities $\rho_j$ are constant we deduce the incompressible flow constraint:
\begin{align}
    {\rm div} \bv_{j}= 0, \quad \text{for }\phi_{j}=1.
\end{align}
The linear momentum of constituent $j$ satisfies the balance law:
\begin{align}\label{eq: lin mom constituent j}
    \partial_t \mathbf{m}_{j} + {\rm div} \left( \mathbf{m}_{j}\otimes \bv_{j} \right) = {\rm div} \mathbf{T}_{j} + \trho_{j} \mathbf{b}_{j} + \bpi_{j} + \bv_{j} \gamma_{j}.
\end{align}
Here $\mathbf{T}_{j}$ is the Cauchy stress tensor of constituent $j$ and $\mathbf{b}_{j}$ is the external body force. The term $\bpi_{j}$ represents the momentum supply of constituent $j$ by the other constituents (see e.g.   \cite{truesdell1960classical}). 
The balance of angular momentum implies that the Cauchy stress tensors of constituents have the form:
\begin{align}
    \mathbf{N}_{j}=\mathbf{T}_{j}-\mathbf{T}_{j}^T,
\end{align}
where $\mathbf{N}_{j}$ represents the intrinsic moment of momentum vector. Note that the last member on the right-hand side of (\ref{eq: lin mom constituent j}) vanishes when switching to a convective form with the aid of the local mass balance \eqref{eq: local mass balance constituent j}:
\begin{align}\label{eq: lin mom constituent j convective}
    \trho_{j}  \grave{\bv}_{j} = {\rm div} \mathbf{T}_{j} + \trho_{j} \mathbf{b}_{j} + \bpi_{j}.
\end{align}

\subsection{Balance laws of mixtures}\label{sec: BL mixtures}
The balance of mass of the mixture density follows by summing \eqref{eq: local mass balance constituent j} over the constituents:
\begin{align}\label{eq: mass balance mixture}
    \partial_t \rho + {\rm div}(\rho \bv) = 0,
\end{align}
where we have postulated the sum of the mass fluxes to vanish:
\begin{align}\label{eq: sum gamma}
    \displaystyle\sum_{j} \gamma_{j} =0.
\end{align}
The mixture velocity $\mathbf{v}$ is in general not divergence-free and this property is in literature referred to as a \textit{quasi-incompressible} mixture \cite{lowengrub1998quasi}.
The postulate \eqref{eq: sum gamma} complies with the third metaphysical principle of mixture theory and precludes the creation of mixture mass. We denote the difference of the mass fluxes by $\gamma = \gamma_1-\gamma_2$, which provides:
\begin{align}
    \gamma_1 = \frac{1}{2}\gamma, \quad\quad \gamma_2= -\frac{1}{2}\gamma.
\end{align}
To proceed, we introduce the order parameters (phase-fields) based on the volume fractions and concentrations of the individual constituents: $\phi=\phi_1-\phi_2 \in [-1,1]$ and $c=c_1-c_2 \in [-1,1]$. By recalling \eqref{eq: sum phi} and \eqref{eq: sum c} we deduce the relations:
\begin{subequations}\label{eq: order parameters}
\begin{alignat}{3}
    \phi_1=&~\frac{1+\phi}{2}, \quad\quad \phi_2=&~\frac{1-\phi}{2},\\
    c_1=&~\frac{1+c}{2}, \quad\quad c_2=&~\frac{1-c}{2}.
\end{alignat}
\end{subequations}
Additionally, the density of the mixture $\rho$ can be expressed in terms $\phi$ and $c$ via $\rho=\hat{\rho}(\phi) = \check{\rho}(c)$ where $\hat{\rho}(\phi)$ and $\check{\rho}(c)$ are defined as:
\begin{subequations}
\begin{alignat}{3}
    \hat{\rho}(\phi) =&~ \rho_1\frac{1+\phi}{2}+\rho_2\frac{1-\phi
    }{2},\\
    \frac{1}{\check{\rho}(c)} =&~ \frac{1}{\rho_1}\frac{1+c}{2} + \frac{1}{\rho_2}\frac{1-c}{2}.
\end{alignat}
\end{subequations}
The relation between $\phi$ and $c$ is given by:
  \begin{align}\label{eq: c phi}
    c = \dfrac{\jrho+\arho \phi}{\arho + \jrho\phi}, \quad 
    \phi = \dfrac{-\jrho+\arho c}{\arho - \jrho c},\quad 
    c'(\phi) = \dfrac{\rho_1\rho_2}{\rho^2}, \quad 
    \phi'(c) = \dfrac{\rho^2}{\rho_1\rho_2}.
  \end{align}
The \textit{diffusive fluxes} are defined as:
\begin{subequations}
\begin{align}
    \bh^v :=&~ \phi_1\bw_1 - \phi_2 \bw_2,\\
    \bJ^v :=&~ \trho_1 \bw_1-\trho_2\bw_2. 
\end{align}
\end{subequations}
The evolution equations of the order parameters $\phi$ and $c$ follow from taking the difference of the mass balance equations of the constituents \eqref{eq: local mass balance constituent j}:
\begin{subequations}\label{eq: phase equations}
\begin{align}
    \dot{\phi} +  \phi {\rm div} \bv + {\rm div} \bh^v=&~ \zeta \gamma, \label{eq: phase equations: phi}\\  \rho\dot{c} +  {\rm div} \bJ^v=&~ \gamma, \label{eq: phase equations: c}
\end{align}
\end{subequations}
where we recall that the dot $\dot{}$ denotes the material derivative with respect to the mixture velocity $\bv$. Next, to obtain the linear momentum equation of the mixture we take the sum of \eqref{eq: lin mom constituent j}:
\begin{align}\label{eq: lin mom mix mass}
    \partial_t \mathbf{m} + {\rm div} \left( \mathbf{m}\otimes \bv \right) = {\rm div} \mathbf{T} + \rho \mathbf{b},
\end{align}
where the Cauchy stress tensor and the body force of the mixture are respectively identified as:
\begin{subequations}
\begin{align}
  \mathbf{T} :=&~ \sum_{j} \mathbf{T}_{j}-\tilde{\rho}_{j}\bw_{j}\otimes\bw_{j},\\
    \rho\mathbf{b} :=&~\sum_{j} \tilde{\rho}_{j}\mathbf{b}_{j}.   
\end{align}
\end{subequations}
Here we have postulated the balance of momentum supplies:
\begin{align}
      \displaystyle\sum_{j} \bpi_{j}+ \gamma_{j}\bv_{j} =0,
\end{align}
which is consistent with the third metaphysical principle of mixture theory and precludes the creation of mixture momentum.
We denote the difference of the growth of linear momentum of the constituents by $\mathbf{p} := (\bpi_1+ \gamma_1\bv_1)-(\bpi_2+ \gamma_2\bv_2)$, which leads to:
\begin{align}
    \bpi_1+ \gamma_1\bv_1 = \frac{1}{2}\mathbf{p} , \quad\quad \bpi_2+ \gamma_2\bv_2= -\frac{1}{2}\mathbf{p} .
\end{align}
Symmetry of the dyadic product implies symmetry of the Cauchy stress tensor of the mixture:
\begin{align}\label{eq: stress sym}
    \mathbf{T}^T = \mathbf{T},
\end{align}
where we have postulated the balance of the intrinsic moment of momentum vectors:
\begin{align}
      \displaystyle\sum_{j} \mathbf{N}_{j} =0.
\end{align}
The corresponding convective form reads:
\begin{align}\label{eq: lin mom mixture convective}
    \rho \dot{\bv}  = {\rm div} \mathbf{T}+\rho\mathbf{b}.
\end{align}

\subsection{Non-typical identities and evolution equations}\label{sec: alternative formulations}
Apart from the mixture velocity \eqref{eq: mix velo} one might define other mean velocities. A typical example of another mean velocity is the \textit{volume-averaged velocity} $\bu$ given by:
\begin{align}\label{eq: volume-averaged velocity}
    \bu(\bx,t) := \displaystyle\sum_{j} \phi_{j}(\bx,t) \bv_{j}(\bx,t).
\end{align}
We define the peculiar velocity relative to the volume-averaged velocity:
\begin{align}\label{eq: def bomegaj}
    \bomega_{j}(\bx,t):=\bv_{j}(\bx,t)-\bu(\bx,t),
\end{align}
which has the consequence:
\begin{align}\label{eq: sum diffusion velocity}
    \displaystyle\sum_{j} \phi_{j}(\bx,t)\bomega_{j}(\bx,t)  = 0.
\end{align}
A direct consequence of the constituent mass balance equations \eqref{eq: local mass balance constituent j} and \eqref{eq: sum diffusion velocity} is the relation:
\begin{align}\label{eq: div free u}
    {\rm div}  \bu = \beta \gamma.
\end{align}
Thus, in absence of mass fluxes ($\gamma = 0$), the volume-averaged velocity inherits the single-constituent solenoidal property of the individual constituents. We additionally introduce the diffusive fluxes with respect to the volume-averaged velocity:
\begin{subequations}\label{eq: diff flux wrt u}
\begin{align}
    \mathbf{h}^u  :=&~ \phi_1\bomega_1-\phi_2\bomega_2,\\
    \bJ^u :=&~ \trho_1 \bomega_1-\trho_2\bomega_2,\\
    \tilde{\bJ}^{u} :=&~ \trho_1\bomega_1+\trho_2\bomega_2.\label{eq: def J tilde}
\end{align}
\end{subequations}
Substitution of the mixture velocity \eqref{eq: mix velo} and the volume-averaged velocity \eqref{eq: volume-averaged velocity} into \eqref{eq: def J tilde} provides the key identity that reveals the difference between the momentum $\mathbf{m}$ and $\rho \bu$:
\begin{align}\label{eq: relation mixture velo, volume avg} 
    \mathbf{m} = \rho \bv = \rho \bu + \tilde{\bJ}^u.
\end{align}
\begin{remark}[Mixture momentum]
  The mixture momentum $\mathbf{m}$ is the sum of the constituent momenta $\mathbf{m}_{j}$ and thus complies with the first metaphysical principle of mixture theory. Relation \eqref{eq: relation mixture velo, volume avg} distinguishes the momentum of the mixture $\mathbf{m}$ from $\rho \bu$. We return to this observation later in this section in the context of volume-averaged velocity NSCH models. \demo
\end{remark}

Substituting the key identity \eqref{eq: relation mixture velo, volume avg} and the property of the velocity field $\bu$ \eqref{eq: div free u} into \eqref{eq: mass balance mixture} provides the alternative form of the mixture mass balance:
\begin{align}
    \overset{\circ}{\rho} + {\rm div} \tilde{\bJ}^u =-\rho\beta\gamma,
\end{align}
where $\circ$ denotes the material derivative with respect to the volume-averaged velocity $\bu$. Furthermore, substituting the volume-averaged velocity $\bu$ and diffusive fluxes \eqref{eq: diff flux wrt u} into \eqref{eq: phase equations} yields the alternative forms of the phase equations::
\begin{subequations}
\begin{align}\ 
    \overset{\circ}{\phi} +  {\rm div} \mathbf{h}^u =&~\dfrac{\rho}{2\rho_1\rho_2}\gamma,\\ 
    \rho\overset{\circ}{c} - c {\rm div} \tilde{\bJ}^u +  {\rm div} \bJ^u=&~ \gamma.
\end{align}
\end{subequations}
With the aim of unifying the various formulation of NSCH models in \cref{sec: 2nd law}, we relate the various diffusive fluxes in the following lemma.
 \begin{lemma}[Relations diffusive fluxes]\label{lem: relations h}
The various diffusive fluxes are related by the identities:
\begin{align}
    \tilde{\bJ}^u = \jrho \bh^u, \quad\quad
    \bJ^u = \left\{\rho\right\}\bh^u,\quad\quad
    \bh^v = \frac{\left\{\rho\right\}}{\rho} \bh^u,\quad\quad
    \bJ^v = \frac{\rho_1\rho_2}{\rho}\bh^u.
\end{align}
\end{lemma}
\begin{proof}
These identities are all consequences of \eqref{eq: sum diffusion velocity 2} and \eqref{eq: sum diffusion velocity}. For example, to obtain the first identity, partition $\tilde{\bJ}^u$ as:
\begin{align}
    \tilde{\bJ}^u =&~ \trho_1\bomega_1+\trho_2\bomega_2\nn\\
    =&~ \frac{\rho_1}{2}(\phi_1\bomega_1-\phi_2\bomega_2)-\frac{\rho_2}{2}(\phi_1\bomega_1-\phi_2\bomega_2)\nn\\
    &~+\frac{\rho_1}{2}(\phi_1\bomega_1+\phi_2\bomega_2)+\frac{\rho_2}{2}(\phi_1\bomega_1+\phi_2\bomega_2).
\end{align}
The first line collapses to $\jrho \bh^u$ and the second line vanishes due to \eqref{eq: sum diffusion velocity}.
\end{proof}
Analogously to the mass balance and phase equations, we may also formulate the momentum balance in terms of the volume-averaged velocity. To this purpose, we first introduce the following lemma.
 \begin{lemma}[Relation peculiar velocities to stress tensor]\label{lem: id stress0}
  We have the identity:
  \begin{align}\label{lem: id stress00}
      \sum_{j} \tilde{\rho}_{j}\bw_{j}\otimes\bw_{j}  =\sum_{j} \tilde{\rho}_{j}\bomega_{j}\otimes\bomega_{j} -\frac{1}{\rho}\tilde{\bJ}^u\otimes\tilde{\bJ}^u.
  \end{align}
\end{lemma}
\begin{proof}
  Identity \eqref{lem: id stress00} is a direct consequence of the following two identities:
  \begin{subequations}
    \begin{align}
      \sum_{j} \tilde{\rho}_{j} \bw_{j} \otimes \bw_{j}  =&~  \sum_{j} \tilde{\rho}_{j} \bv_{j} \otimes \bv_{j} - \rho \bv \otimes \bv, \label{eq: id1}\\ 
      \sum_{j} \tilde{\rho}_{j} \bomega_{j} \otimes \bomega_{j} =&~\sum_{j} \tilde{\rho}_{j} \bv_{j} \otimes \bv_{j} - \rho \bv \otimes \bv + \frac{1}{\rho}\tilde{\bJ}^u\otimes\tilde{\bJ}^u. \label{eq: id2}
    \end{align}
  \end{subequations}
  To see the first identity, we substitute the definition of the diffuse-velocity $\bw_{j}$ \eqref{eq: def bwj} into the left-hand side of \eqref{eq: id1} and expand:
  \begin{align}
    \sum_{j} \tilde{\rho}_{j} \bw_{j} \otimes \bw_{j} =&~ \sum_{j} \tilde{\rho}_{j} (\bv_{j}-\bv) \otimes (\bv_{j}-\bv)\nn\\
    =&~ \sum_{j} \left(\tilde{\rho}_{j} \bv_{j} \otimes \bv_{j} - \tilde{\rho}_{j} \bv_{j} \otimes \bv - \tilde{\rho}_{j} \bv \otimes \bv_{j} +\tilde{\rho}_{j} \bv \otimes \bv \right) \nn\\   
    =&~ \sum_{j} \tilde{\rho}_{j} \bv_{j} \otimes \bv_{j} - \rho \bv \otimes \bv - \rho \bv \otimes \bv +\rho \bv \otimes \bv  \nn\\
    =&~ \sum_{j} \tilde{\rho}_{j} \bv_{j} \otimes \bv_{j} - \rho \bv \otimes \bv.
  \end{align}
  The second follows in the same fashion via the substitution of the definition of $\bomega_{j}$ \eqref{eq: def bomegaj} and the key identity \eqref{eq: relation mixture velo, volume avg} into the left-hand side of \eqref{eq: id2}:  
  
  \begin{align}
    \sum_{j} \tilde{\rho}_{j} \bomega_{j} \otimes \bomega_{j} =&~ \sum_{j} \tilde{\rho}_{j} (\bv_{j}-\bu) \otimes (\bv_{j}-\bu)\nn\\
    =&~ \sum_{j} \left(\tilde{\rho}_{j} \bv_{j} \otimes \bv_{j} - \tilde{\rho}_{j} \bv_{j} \otimes \bu - \tilde{\rho}_{j} \bu \otimes \bv_{j} +\tilde{\rho}_{j} \bu \otimes \bu \right) \nn\\
    =&~ \sum_{j} \tilde{\rho}_{j} \bv_{j} \otimes \bv_{j} - \rho \bv \otimes \bu - \rho \bu \otimes \bv +\rho \bu \otimes \bu  \nn\\
    =&~ \sum_{j} \tilde{\rho}_{j} \bv_{j} \otimes \bv_{j} - \rho \bv \otimes \bv + \frac{1}{\rho}\tilde{\bJ}^u\otimes\tilde{\bJ}^u.
  \end{align}\\
\end{proof}

By substituting the key identity \eqref{eq: relation mixture velo, volume avg} into the mixture linear momentum equation \eqref{eq: lin mom mix mass} and making use of \cref{lem: id stress0}, one can show that the mixture linear momentum equation formulated in terms of the volume-averaged velocity takes the following form:
\begin{align}\label{eq: mix mom vol avg velo}
  \partial_t (\rho \bu + \tilde{\bJ}^u) + {\rm div} \left( \rho \bu\otimes \bu + \tilde{\bJ}^u\otimes \bu + \bu\otimes \tilde{\bJ}^u + \frac{1}{\rho}\tilde{\bJ}^u\otimes\tilde{\bJ}^u\right) = {\rm div} \mathbf{T}+\rho\mathbf{g}.
\end{align}

Analogously to taking the difference of the constituent mass conservation equations, we may evaluate the difference of the constituent linear momentum equations. This provides a balance law for the diffusive flux. To do so, we first note that by \eqref{eq: local mass balance constituent j conv} and \eqref{eq: lin mom constituent j convective} we can rewrite the constituent balance law of linear momentum as:
\begin{align}\label{eq: rewrite BL of linear momentum}
    \left(\trho_{j}\bv_{j}\right)\grave{} =  \divg \mathbf{T}_{j} + \trho_j\mathbf{b}_{j} + \bpi_{j}+ \gamma_{j}\bv_{j} - \trho_j \bv_j \divg \bv_j.
\end{align}
Denoting $\bJ_j^v=\trho_{j}\bw_{j}$, we have the identity:
\begin{align}
  \left(\trho_{j}\bv_{j}\right)\grave{} =  \grave{\bJ}_j^v+\grave{\trho}_{j}\mathbf{v} +(\nabla \bv) \bJ_j^v + \trho_{j}\dot{\bv}.
\end{align}
Substituting this into \eqref{eq: rewrite BL of linear momentum} provides: 
\begin{align}
    \grave{\bJ}_j^v + (\nabla \bv) \bJ_j^v=  \divg \mathbf{T}_{j} + \trho_j\mathbf{b}_{j} + \bpi_{j}+ \gamma_{j}\bv_{j} -\grave{\trho}_{j}\mathbf{v}   - \trho_{j}\dot{\bv}- \trho_j \bv_j \divg \bv_j. 
\end{align}
By expanding the peculiar derivative as:
  \begin{align}
    \grave{\bJ}_j^v =&~\dot{\bJ}_j^v + \bw_{j}\cdot\nabla \bJ_j^v\nn\\
    =&~\dot{\bJ}_j^v + \divg\left(\bJ_j^v\otimes \bw_j\right) -\left(\divg \bv_{j}\right)\bJ_j^v +(\divg \bv) \bJ_j^v,
  \end{align}
we arrive at:
\begin{align}\label{eq: diff flux intermediate}
    \dot{\bJ}_j^v + \left(\nabla \bv +(\divg \bv)\mathbf{I}\right) \bJ_j^v =&~  \divg \left(\mathbf{T}_{j}-\bJ_j^v\otimes \bw_j\right) + \trho_j\mathbf{b}_{j} \nonumber \\
    &~+ \bpi_{j}+ \gamma_{j}\bv_{j} -\gamma_j\bv- \trho_{j}\dot{\bv},
\end{align}
where we have used the identity:
\begin{align}
  -\grave{\trho}_{j}\mathbf{v}  +\left(\divg \bv_{j}\right)\bJ_j^v- \trho_j \bv_j \divg \bv_j = -\gamma_j\bv.
\end{align}
By substituting for the last member in \eqref{eq: diff flux intermediate} the mixture linear momentum evolution equation \eqref{eq: lin mom mixture convective} in the form
\begin{align}
        \trho_{j} \dot{\bv}  = c_j{\rm div} \mathbf{T}+\trho_{j}\mathbf{b},
\end{align}
we obtain:
\begin{align}\label{eq: evo Jv}
    \dot{\bJ}^v_j + \left(\nabla \bv +(\divg \bv)\mathbf{I}\right) \bJ_j^v=&~ \divg \left(\mathbf{T}_{j}-\bJ_j^v\otimes \bw_j\right) - c_j \divg \mathbf{T} \nn\\
    &~ + \bpi_j+ \gamma_j \bv_j-\gamma_j\bv+ \trho_j \mathbf{b}_j -\trho_{j} \mathbf{b}.
\end{align}
Note that in case the body force is the same for all components, i.e. $\mathbf{b}_{j}=\mathbf{b}$, the body force drops out of the evolution equation \eqref{eq: evo Jv}.
By subtracting \eqref{eq: evo Jv} for constituent $2$ from that of constituent $1$ we arrive at \textit{the evolution equation of the diffusive flux}:
\begin{align}\label{eq: diff flux eq}
    \dot{\bJ}^v + \left( \nabla \bv + (\divg \bv)\mathbf{I}  \right) \bJ^v=&~ \divg \left(\mathbf{T}_1-\mathbf{T}_2 + \dfrac{\tilde{\rho}_1-\tilde{\rho}_2}{4\trho_1\trho_2}\bJ^v\otimes \bJ^v\right) - c \divg \mathbf{T} \nn\\
    &~+ \trho_1 (\mathbf{b}_1-\mathbf{b})-\trho_2 (\mathbf{b}_2-\mathbf{b})  + \mathbf{p}-\gamma \mathbf{v},
\end{align}
where we have utilized the identity:
\begin{align}
  -\bJ_1^v\otimes \bw_1+\bJ_2^v\otimes \bw_2 = \dfrac{\tilde{\rho}_1-\tilde{\rho}_2}{4\trho_1\trho_2}\bJ^v\otimes \bJ^v, \quad \text{ for }\phi_1\phi_2\neq 0,
\end{align}
that results from \eqref{eq: sum diffusion velocity 2}. This term vanishes for $\phi_1\phi_2= 0$.

With the aid of \cref{lem: relations h}, the key identity \eqref{eq: relation mixture velo, volume avg} and the mixture mass balance \eqref{eq: mass balance mixture}, one can formulate the evolution equation of the diffusive flux in terms of $\bJ^u$. First, the material derivative of the diffusive flux $\bJ^v$ may be expressed in terms of $\bJ^u$ via:
\begin{align}\label{eq: diff flux 6}
    \dot{\bJ}^v = \zeta^{-1}\rho^{-1}\left(\dot{\bJ}^u-\alpha\bJ^u\divg \left( \rho^{-1}\bJ^u\right)+\beta\gamma\bJ^u\right).
\end{align}
Next, we may express the material derivative $\dot{\bJ}^u$ as:
\begin{align}
    \dot{\bJ}^u=\overset{\circ}{\bJ^u}-\alpha\rho^{-1}\bJ^u\cdot\nabla \bJ^u,
\end{align}
and substitution into \eqref{eq: diff flux 6} provides:
\begin{align}\label{eq: diff flux 7}
    \dot{\bJ}^v = \zeta^{-1}\rho^{-1}\left(\overset{\circ}{\bJ^u}-\alpha\divg \left( \rho^{-1}\bJ^u\otimes\bJ^u \right)+\beta\gamma\bJ^u\right).
\end{align}
Thus the evolution equation of $\bJ^u$ may be written as:
\begin{align}\label{eq: diff flux eq u}
    &\overset{\circ}{\bJ^u}-\alpha\divg \left( \rho^{-1}\bJ^u\otimes\bJ^u \right) +\beta\gamma\bJ^u+\left(\nabla\left(\bu-\alpha\rho^{-1}\bJ^u\right)+\divg\left(-\alpha\rho^{-1}\bJ^u\right)\right)\bJ^u=\nn\\ &\quad \zeta\rho\left(\divg \left(\mathbf{T}_1-\mathbf{T}_2+\zeta^{-2}\rho^{-2}\dfrac{\tilde{\rho}_1-\tilde{\rho}_2}{4\trho_1\trho_2}\bJ^u\otimes \bJ^u\right) - c \divg \mathbf{T} \right.\nonumber\\
    &\quad\quad\quad\quad\left.+ \trho_1 (\mathbf{b}_1-\mathbf{b})-\trho_2 (\mathbf{b}_2-\mathbf{b})  + \mathbf{p}-\gamma \mathbf{v}\right).
\end{align}
\begin{remark}[Mixture momentum equation]
  An equivalent form of \eqref{eq: mix mom vol avg velo} has appeared in the Ph.D. thesis of Simsek \cite{simsek2017thesis}. To the best knowledge of the authors, the linear mixture momentum equation of NSCH models formulated in terms of the volume-averaged velocity appears in two distinct forms. The first is of the form \eqref{eq: lin mom mix mass} with $\bu$ instead of the mixture velocity $\bv$ and $\rho \bu$ instead of $\mathbf{m}$. The other, proposed by Abels et al. \cite{abels2012thermodynamically}, follows, as observed by Simsek \cite{simsek2017thesis}, when taking $\overset{\circ}{\bJ^u} = 0$ in \eqref{eq: mix mom vol avg velo}. From the standpoint of mixture theory, the first class relies on the hidden assumption that all the terms containing $\tilde{\bJ}^u$ vanish, whereas the model of Abels et al. contains the hidden assumption $\overset{\circ}{\bJ^u} = 0$. Both assumptions are incompatible with \eqref{eq: diff flux eq} and the linear mixture momentum equation in these models does not match with mixture theory.
  \demo
\end{remark}
\begin{remark}[Evolution equation diffusive flux]
  The evolution equation of the diffusive flux seems to be not well known in the phase-field community. We note however that it was first derived in 1975 by M\"{u}ller \cite{muller1975thermo} and revisited in e.g.  \cite{morro2016nonlinear} and   \cite{muller2001thermodynamics}. \demo
\end{remark}

To summarize, in agreement with the first unification principle, we have obtained the following \textit{equivalent formulations of the same model}, one formulated in terms of the mixture velocity $\bv$:
\begin{subequations}\label{eq: model 00}
  \begin{empheq}[left=\empheqlbrace]{align}
    \partial_t (\rho \bv) + {\rm div} \left( \rho \bv\otimes \bv \right) - {\rm div} \mathbf{T}-\rho\mathbf{g}&=~ 0, \label{eq: model 00: mom}\\
 \partial_t \rho + {\rm div}(\rho \bv) &=~ 0, \label{eq: model 00: cont} \\
 \dot{\bJ}^v + \left( \nabla \bv + (\divg \bv)\mathbf{I}  \right) \bJ^v- \divg \left(\mathbf{T}_1-\mathbf{T}_2 + \dfrac{\tilde{\rho}_1-\tilde{\rho}_2}{4\trho_1\trho_2}\bJ^v\otimes \bJ^v\right)  &\nn\\
 + c \divg \mathbf{T}- \trho_1 (\mathbf{b}_1-\mathbf{b})+\trho_2 (\mathbf{b}_2-\mathbf{b})  -( \mathbf{p}-\gamma \mathbf{v})&=~ 0, \label{eq: model 00: h} \\
 \text{with the evolution equation of an order parameter,}&\nn\\
 \text{either:} \quad \dot{\phi} +  \phi {\rm div} \bv + {\rm div} \bh^v -\zeta\gamma&=~0,  \label{eq: model 00: PF}\\
 \text{or: }\quad\quad\quad\quad\quad \rho\dot{c} +  {\rm div} \bJ^v
  - \gamma&=~0,
  \end{empheq}
\end{subequations}
and one formulated in terms of the volume-averaged velocity $\bu$:
\begin{subequations}\label{eq: model 0}
  \begin{empheq}[left=\empheqlbrace]{align}
    \partial_t (\rho \bu + \tilde{\bJ}^u) + {\rm div} \left( \rho \bu\otimes \bu + \tilde{\bJ}^u\otimes \bu + \bu\otimes \tilde{\bJ}^u + \frac{1}{\rho}\tilde{\bJ}^u\otimes\tilde{\bJ}^u \right)&\nn\\
    - {\rm div} \mathbf{T}-\rho\mathbf{g}&=~ 0, \label{eq: model 0: mom}\\
 {\rm div} \bu-\beta\gamma &=~ 0, \label{eq: model 0: cont} \\
 \overset{\circ}{\bJ^u}-\alpha\divg \left( \rho^{-1}\bJ^u\otimes\bJ^u \right)+\beta\gamma\bJ^u&\nn\\
 +\left(\nabla\left(\bu-\alpha\rho^{-1}\bJ^u\right)+\divg\left(-\alpha\rho^{-1}\bJ^u\right)\right)\bJ^u&\nn\\ -\frac{\beta\rho}{\alpha}\left(\divg \left(\mathbf{T}_1-\mathbf{T}_2+\alpha^2\beta^{-2}\rho^{-2}\dfrac{\tilde{\rho}_1-\tilde{\rho}_2}{4\trho_1\trho_2}\bJ^u\otimes \bJ^u\right) - c \divg \mathbf{T} \right.&\nn\\
 \left.+\trho_1 (\mathbf{b}_1-\mathbf{b})-\trho_2 (\mathbf{b}_2-\mathbf{b})   + \mathbf{p}-\gamma \mathbf{v}\right) &=~0,\\
 \text{with the evolution equation of an order parameter,}&\nn\\
 \text{either:} \quad\quad\quad\quad\quad \overset{\circ}{\phi}+ {\rm div} \bh^u- \dfrac{\rho}{2\rho_1\rho_2}\gamma&=~ 0,\label{eq: model 0: PF}\\
 \text{or:} \quad \rho\overset{\circ}{c} - c {\rm div} \tilde{\bJ}^u +  {\rm div} \bJ^u -\gamma
 &=~ 0.
  \end{empheq}
\end{subequations}

\section{Constitutive modeling}\label{sec: 2nd law}
In this section we discuss the constitutive modeling. We will use the Coleman-Noll procedure \cite{coleman1974thermodynamics} to obtain constitutive equations that guarantee an energy-dissipation equation.
\cref{sec: dam} provides the necessary definitions and assumptions. In \cref{sec: Derivation of the constitutive modeling restriction} we establish the constitutive modeling restriction. Then, in \cref{subsec: Modeling restrictions} we discuss the equivalence of alternative modeling restrictions. The actual selection of the constitutive models appears in \cref{sec: Selection of constitutive models}.

\subsection{Definitions, assumptions and modeling choices}\label{sec: dam}

We make the following \textit{simplification assumptions} (S):
\begin{itemize}
    \item the kinetic energies of the constituents are negligible when computed relative to the gross motion of the constituent.
    \item the body force acting on the constituents is a constant gravitational force, i.e. $\mathbf{b}_j=\mathbf{b} = \bg = -g \boldsymbol{\jmath}$ with $g$ constant and $\boldsymbol{\jmath}$ the vertical unit vector.
\end{itemize} 
We define the total energy $\mathscr{E}$ associated with the system of balance laws as the sum of the Helmholtz free energy, kinetic energy and gravitational energy:
\begin{align}\label{eq: total energy}
  \mathscr{E}(\phi,\bv) =  \displaystyle\int_{\mathcal{R}(t)}(\Psi + \mathscr{K}(\mathcal{R}(t)) + \mathscr{G}(\mathcal{R}(t)))~{\rm d}v.
\end{align}
Here $\mathcal{R}=\mathcal{R}(t)$ denotes an arbitrary time-dependent control volume in $\Omega$ with volume element ${\rm d}v$ and unit outward normal $\boldsymbol{\nu}$ that is transported by $\bv$ (and thus the normal velocity of $\partial \mathcal{R}(t)$ is $\bv\cdot \boldsymbol{\nu}$). The kinetic and gravitational energies are given by:
\begin{subequations}
\begin{align}
    \mathscr{K}(\mathcal{R}(t)) =&~ \frac{1}{2} \rho \|\bv\|^2, \label{eq: def Ekin}\\
    \mathscr{G}(\mathcal{R}(t)) =&~ \rho g y.
\end{align}
\end{subequations}
\begin{remark}[Free energy]
  The free energy $\Psi$ is defined with respect to the volume element ${\rm d}v$. We call this free energy a \textit{volume-measure-based} free energy. The common alternative is to introduce a free energy $\psi$ with respect to the mass element $\rho{\rm d}v$. This free energy is referred to as a \textit{mass-measure-based} free energy. Throughout this paper we use $\Psi$ and $\psi$ for volume-measure-based, and mass-measure-based free energies, respectively.
\end{remark}
\begin{remark}[Deviation from mixture theory]\label{rmk: Deviation from mixture theory}
At this point we deviate from mixture theory of rational mechanics laid down by Truesdell \cite{truesdell1984historical}. Specifically, we diverge from mixture theory in that the first metaphysical principle is violated. Namely, mixture theory dictates that the kinetic energy of the mixture $\mathscr{K}_M$ is the superposition of the individual kinetic energies of the constituents:
  \begin{align}
      \mathscr{K}_M = \displaystyle\sum_{j}\mathscr{K}_{j}.
  \end{align}
  A straightforward substitution of the mixture quantities $\rho$ and $\rho \bv$ reveals:
  \begin{align}\label{eq: link KM and K}
      \mathscr{K}_M = \mathscr{K} + \displaystyle\sum_{j} \frac{1}{2}\trho_{j} \|\bw_{j}\|^2,
  \end{align}
  in which the second member represents the kinetic energy of the constituents relative to the gross motion of the constituent. Remark that the second member vanishes in the single constituent and thus the discrepancy appears solely in the interface region. The simplification assumption (S) requires the vanishing of the second member. \demo
\end{remark}

\begin{remark}[Simplification assumption]
The first part of the simplification assumption (S) is closely related to the one adopted by Gurtin \cite{gurtinmodel} which states: \textit{`the momenta and kinetic energies of the constituents are negligible when computed relative to the gross motion of the constituent'}. We recall that the first part is not an assumption since the momenta relative to the gross motion of the constituent is absent as stated in \eqref{eq: rel gross motion zero}.
\demo
\end{remark}

\begin{remark}[Kinetic energy volume-averaged velocity models]\label{rmk: Ekin}
  The kinetic energy of diffuse-interface models formulated in terms of the volume-averaged velocity $\bu$ is taken as $\rho \|\bu\|^2/2$, see e.g.  Abels et al. \cite{abels2012thermodynamically}. Just as using $\mathscr{K}$ instead of the kinetic energy of the mixture $\mathscr{K}_M$, see \cref{rmk: Deviation from mixture theory}, the mismatch exclusively occurs in the interface region. An important difference is that a relation like \eqref{eq: link KM and K}, in which the kinetic energy of the mixture is the sum of the approximate kinetic energy and a kinetic energy of with respect to the gross motion of the constituent, does not hold. Instead we have the relation:
    \begin{align}\label{eq: link KM and K(u)}
      \mathscr{K}_M = \frac{1}{2} \rho \|\bu\|^2 + \displaystyle\sum_{j} \frac{1}{2}\trho_{j} \|\bw_{j}\|^2  + \frac{1}{2\rho}\tilde{\bJ}^u\cdot\tilde{\bJ}^u+ \tilde{\bJ}^u \cdot \bu,
  \end{align}
  which follows from the identity relation \eqref{eq: relation mixture velo, volume avg} and reveals a non obvious approximation of the kinetic energy of the mixture by $\rho \|\bu\|^2/2$. The last term in \eqref{eq: link KM and K(u)} is not guaranteed positive and thus does not represent a kinetic energy relative to some different velocity. Relation \eqref{eq: relation mixture velo, volume avg} reveals that the kinetic energy $\mathscr{K}$ formulated in terms of the volume-averaged velocity reads:
  \begin{align}
    \mathscr{K}(\mathcal{R}(t)) =&~ \dfrac{(\rho \bu + \tilde{\bJ}^u)\cdot(\rho \bu + \tilde{\bJ}^u) }{2 \rho}.
\end{align}\demo
\end{remark}
We proceed with selecting the volume fraction difference $\phi$ as order parameter and postulating the free energy to pertain to the constitutive class:
\begin{align}\label{eq: class Psi}
  \Psi = \hat{\Psi}(\phi,\nabla \phi,\mathbf{D}),   
\end{align}
where we recall that $\mathbf{D}$ is the symmetric gradient of the mixture velocity:
\begin{align}
  \mathbf{D} = \frac{1}{2}\left( \nabla \bv + (\nabla \bv)^T \right).
\end{align}
Next, we define a chemical potential-like quantity in the usual way as:
\begin{align}
    \hat{\mu} := \dfrac{ \partial \hat{\Psi}}{\partial \phi} - {\rm div}\dfrac{\partial \hat{\Psi}}{\partial\nabla \phi}.
\end{align}
In the phase-field community this is sometimes referred to as the Fr\'{e}chet (or variational) derivative of the total Helmholtz free energy. We postpone the selection of the constitutive classes of the stress tensor, diffusive flux and the mass flux, denoted $\hat{\mathbf{T}}, \hat{\bh}^v$ and $\hat{\gamma}$, to the end of \cref{sec: Derivation of the constitutive modeling restriction}.
With the aim of deriving NSCH models, we postulate the energy-dissipation law:
\begin{align}\label{eq: second law}
    \dfrac{{\rm d}}{{\rm d}t} \mathscr{E}(\phi,\bv) = \mathscr{W}(\mathcal{R}(t)) - \mathscr{D}(\mathcal{R}(t)),
\end{align}
where $\mathscr{W}(\mathcal{R}(t))$ is the rate of work performed by macro- and micro stresses coming through the boundary $\partial \mathcal{R}(t)$ and $\mathscr{D}(\mathcal{R}(t))$ is the dissipation for which we demand $\mathscr{D}(\mathcal{R}(t))\geq 0$.

\begin{remark}[Arbitrariness of modeling choices]\label{rmk: arbitrariness of modeling choices}
The choice of working with the volume fraction difference $\phi$ as the order parameter seems arbitrary. One could instead work with the concentration difference $c$. Additionally, one may work with a mass-measure-based free energy instead of a volume-measure-based free energy. In \cref{subsec: Modeling restrictions} we discuss the relation of the modeling choices in detail.
\demo
\end{remark}
\begin{remark}[Constitutive class diffusive flux]
From the viewpoint of mixture theory, the usage of a constitutive class for the diffusive flux should come as a surprise. Namely, in \cref{sec: mix theory} we have established an evolution equation for the mass diffusion $\bJ^v$, i.e. \eqref{eq: diff flux eq}. Here we discard the PDE \eqref{eq: diff flux eq}. This is an approximation of the mixture theory framework, or at least, it is not obvious how it fits in the framework. This approximation however is essential in order to work with a reduced model of NSCH type. In other words, without the approximation one does not retrieve a NSCH type model.
\demo
\end{remark}
\begin{remark}[Energy-dissipation law]
The energy-dissipation statement is a necessary and core element in the derivation of NSCH models. The equation \eqref{eq: second law} is closely linked to the first law of thermodynamics whereas the requirement $\mathscr{D}(\mathcal{R}(t))\geq 0$ is associated with the second law of thermodynamics. In the literature on NSCH models the energy-dissipation law is often called the second law of thermodynamics and models that satisfy it are referred to as thermodynamically consistent. However, we note that the energy-dissipation statement is not obviously compatible with the second law of thermodynamics presented in mixture theory, see e.g.  \cite{truesdell1984historical}. As such referring to NSCH models as thermodynamically consistent is not justified. We refer to  Aki et al. \cite{aki2014quasi} and  \cite{dreyer2010van} for some remarks on approximations in the energy-dissipation statement.
\demo
\end{remark}
\begin{remark}[Rate of work term]
By using concepts like the microforce balance, the rate of work term may be already defined at this point, see  Gurtin \cite{gurtinmodel}.
\demo
\end{remark}
\begin{remark}[Alternative modeling approach]
Instead of working with constitutive classes and the Coleman-Noll procedure, an alternative is to postulate the dissipation, see  \cite{malek2008thermodynamic}.
\demo
\end{remark}
\subsection{Derivation of the constitutive modeling restriction}\label{sec: Derivation of the constitutive modeling restriction}
We proceed with the evaluation of the evolution of the energy \eqref{eq: total energy}. By applying Reynolds transport theorem to the free energy $\hat{\Psi}$ we have:
\begin{align}
      \dfrac{{\rm d}}{{\rm d}t}\displaystyle\int_{\mathcal{R}(t)} \hat{\Psi} ~{\rm d}v = \displaystyle\int_{\mathcal{R}(t)} \partial_t \hat{\Psi} ~{\rm d}v + \displaystyle\int_{\partial \mathcal{R}(t)} \hat{\Psi} \bv \cdot \bnu  ~{\rm d}a.
\end{align}
In the next step we apply the divergence theorem and expand the derivatives:
\begin{align}\label{eq: Psi derivation 1}
    \dfrac{{\rm d}}{{\rm d}t}\displaystyle\int_{\mathcal{R}(t)} \hat{\Psi} ~{\rm d}v  &~= \displaystyle\int_{\mathcal{R}(t)} \dfrac{\partial \hat{\Psi}}{\partial \phi} \dot{\phi} + \dfrac{\partial \hat{\Psi}}{\partial \nabla \phi}\cdot \left(\nabla \phi\right)\dot{}+ \dfrac{\partial \hat{\Psi}}{\partial \mathbf{D}} : \dot{\mathbf{D}} + \hat{\Psi} ~{\rm div} \bv~{\rm d}v,
\end{align}
where we recall that the dot $\dot{}$ is the material derivative with respect to the mixture velocity $\bv$. By substituting the identity
\begin{align}\label{eq: relation grad phi}
    \left(\nabla \phi\right)\dot{} = \nabla (\dot{\phi}) - (\nabla \phi)^T\nabla \bv,
\end{align}
into \eqref{eq: Psi derivation 1} and subsequently integrating by parts, we arrive at:
\begin{align}
    \dfrac{{\rm d}}{{\rm d}t}\displaystyle\int_{\mathcal{R}(t)} \hat{\Psi} ~{\rm d}v  =&~ \displaystyle\int_{\mathcal{R}(t)} \hat{\mu} \dot{\phi} - \left(\nabla \phi \otimes \dfrac{\partial \hat{\Psi}}{\partial \nabla \phi}\right): \nabla \bv + \dfrac{\partial \hat{\Psi}}{\partial \mathbf{D}} : \dot{\mathbf{D}}  + \hat{\Psi}~{\rm div} \bv~{\rm d}v  \nn\\
    &~+ \displaystyle\int_{\partial \mathcal{R}(t)}\dot{\phi} \dfrac{\partial \hat{\Psi}}{\partial \nabla \phi}\cdot \boldsymbol{\nu} ~{\rm d}a.
\end{align}
We substitute the phase evolution equation \eqref{eq: phase equations: phi} for the material derivative $\dot{\phi}$ and apply integration by parts to obtain:
\begin{align}\label{eq: Psi 2}
    \dfrac{{\rm d}}{{\rm d}t}\displaystyle\int_{\mathcal{R}(t)} \hat{\Psi} ~{\rm d}v = &~ \displaystyle\int_{\mathcal{R}(t)} \nabla \hat{\mu} \cdot \hat{\bh}^v  - \left(\nabla \phi \otimes \dfrac{\partial \hat{\Psi}}{\partial \nabla \phi}\right): \nabla \bv + \dfrac{\partial \hat{\Psi}}{\partial \mathbf{D}} : \dot{\mathbf{D}}   \nn\\
    &~\quad \quad + (\hat{\Psi}-\hat{\mu}\phi)~{\rm div} \bv + \hat{\mu} \zeta \hat{\gamma}~{\rm d}v\nn\\
    &~+ \displaystyle\int_{\partial \mathcal{R}(t)}\left(-\hat{\mu} \hat{\bh}^v+\dot{\phi} \dfrac{\partial \hat{\Psi}}{\partial \nabla \phi}\right)\cdot \boldsymbol{\nu} ~{\rm d}a.
\end{align}
Next, we turn our focus on the kinetic energy. In a similar fashion we apply Reynolds transport theorem and find:
\begin{align}
    \dfrac{{\rm d}}{{\rm d}t}\displaystyle\int_{\mathcal{R}(t)} \mathscr{K}(\mathcal{R}(t)) ~{\rm d}v =&~ \displaystyle\int_{\mathcal{R}(t)} \bv\cdot \partial_t (\rho \bv)-\frac{1}{2} \|\bv\|^2\partial_t \rho  ~{\rm d}v \nn\\
    &~+ \displaystyle\int_{\partial \mathcal{R}(t)} \mathscr{K}(\mathcal{R}(t)) \bv \cdot \bnu  ~{\rm d}a.
\end{align}
Substitution of the mass and momentum equations provides:
\begin{align}
    \dfrac{{\rm d}}{{\rm d}t}\displaystyle\int_{\mathcal{R}(t)} \mathscr{K}(\mathcal{R}(t)) ~{\rm d}v =&~ \displaystyle\int_{\mathcal{R}(t)} -\bv\cdot {\rm div} \left( \rho \bv\otimes \bv \right)+\frac{1}{2} \|\bv\|^2{\rm div}~(\rho \bv) + \bv \cdot {\rm div}~  \hat{\mathbf{T}}\nn\\
    &~\quad \quad +\rho\bv\cdot\mathbf{g}    ~{\rm d}v + \displaystyle\int_{\partial \mathcal{R}(t)} \mathscr{K}(\mathcal{R}(t)) \bv \cdot \bnu  ~{\rm d}a.
\end{align}
With the aid of the identity
\begin{align}
    -\bv \cdot {\rm div} \left( \rho \bv \otimes \bv \right)+ \tfrac{1}{2}\|\bv\|_2^2 {\rm div}(\rho\bv)=&~ - {\rm div} \left( \tfrac{1}{2}\rho\|\bv\|^2 \bu\right),
\end{align}
and using integration by parts on the stress tensor, the evolution of the kinetic energy simplifies to:
\begin{align}\label{eq: kin2}
    \dfrac{{\rm d}}{{\rm d}t}\displaystyle\int_{\mathcal{R}(t)} \mathscr{K}(\mathcal{R}(t)) ~{\rm d}v =&~ \displaystyle\int_{\mathcal{R}(t)} - \nabla \bv : \hat{\mathbf{T}}+\rho\bv\cdot\mathbf{g}    ~{\rm d}v+ \displaystyle\int_{\partial \mathcal{R}(t)} \bv \cdot \hat{\mathbf{T}} \bnu  ~{\rm d}a.
\end{align}
Finally, for the gravitational energy evolution we also apply Reynolds transport theorem to find:
\begin{align}\label{eq: grav}
    \dfrac{{\rm d}}{{\rm d}t}\displaystyle\int_{\mathcal{R}(t)} \mathscr{G}(\mathcal{R}(t)) ~{\rm d}v =&~ \displaystyle\int_{\mathcal{R}(t)} g y \partial_t \rho~{\rm d}v + \displaystyle\int_{\partial \mathcal{R}(t)} \mathscr{G}(\mathcal{R}(t)) \bv \cdot \bnu  ~{\rm d}a.
  \end{align}
Substituting the mass evolution equation \eqref{eq: mass balance mixture} and subsequently providing integration by part leads to:
  \begin{align}\label{eq: grav2}
    \dfrac{{\rm d}}{{\rm d}t}\displaystyle\int_{\mathcal{R}(t)} \mathscr{G}(\mathcal{R}(t)) ~{\rm d}v =&~- \displaystyle\int_{\mathcal{R}(t)} \rho \bv\cdot\mathbf{g}~{\rm d}v.
  \end{align}
Taking the sum of \eqref{eq: Psi 2}, \eqref{eq: kin2} and \eqref{eq: grav2} we arrive at:
\begin{align}\label{eq: second law subst 1}
    \dfrac{{\rm d}}{{\rm d}t} \mathscr{E}(\phi,\bv) = &~ \displaystyle\int_{\partial \mathcal{R}(t)}\left(\bv^T\hat{\mathbf{T}}-\hat{\mu} \hat{\bh}^v+\dot{\phi} \dfrac{\partial \hat{\Psi}}{\partial \nabla \phi}\right)\cdot \boldsymbol{\nu} ~{\rm d}a \nn\\
    &~- \displaystyle\int_{\mathcal{R}(t)}  \left( \hat{\mathbf{T}}  + \nabla \phi \otimes \dfrac{\partial \hat{\Psi}}{\partial \nabla \phi} + (\hat{\mu}\phi-\hat{\Psi})\mathbf{I} \right): \nabla \bv -\nabla \hat{\mu} \cdot \hat{\bh}^v\nn\\
    &~\quad \quad \quad - \dfrac{\partial \hat{\Psi}}{\partial \mathbf{D}} : \dot{\mathbf{D}} -\hat{\mu} \zeta \hat{\gamma} ~{\rm d}v.
\end{align}
Next, we observe that the components $\nabla \bv$, $\hat{\bh}^v$ and $\hat{\gamma}$ are not independent. Namely, we have the sequence of identities:
\begin{align}\label{eq: id p0}
    -  ~{\rm div}\bv = -  \alpha \left(  \dot{\phi} + \phi ~{\rm div}\bv \right)=   \alpha ~{\rm div}\hat{\bh}^v- \beta \hat{\gamma}.
\end{align}
The first identity is a direct consequence of the mass balance of the mixture \eqref{eq: mass balance mixture} and the second follows from the phase equation \eqref{eq: phase equations: phi}. We now introduce the following partition:
\begin{align}\label{eq: partition}
    \hat{\mathbf{T}} = \hat{\mathbf{T}}_0 - \hat{p} \mathbf{I} ~\text{with}~\hat{\mathbf{T}}_0 := \hat{\mathbf{T}} + \hat{p} \mathbf{I},
\end{align}
in which $\hat{p}$ is a scalar field that represents the mechanical pressure and $-\hat{p}\mathbf{I}$ corresponds to the hydro-static part of the stress tensor $\hat{\mathbf{T}}$. To exploit the degeneracy, we let the scalar field $\hat{p}$ act as a Lagrange multiplier on \eqref{eq: id p0} by multiplying this equation by $\hat{p}$: \begin{align}\label{eq: id p}
    - \hat{p} ~{\rm div}\bv =&~  \hat{p} \alpha ~{\rm div}\hat{\bh}^v-\hat{p} \beta \hat{\gamma} \nn\\ =&~ -\nabla (\hat{p} \alpha)\cdot \hat{\bh}^v+ {\rm div}(\hat{p} \alpha \hat{\bh}^v)-\hat{p}  \beta \hat{\gamma}.
\end{align}
If we substitute the partition \eqref{eq: partition} and the identity \eqref{eq: id p} into \eqref{eq: second law subst 1} we can identify the rate of work performed by macro- and micro stresses and the dissipation as:
\begin{subequations}\label{eq: W, D}
\begin{align}
    \mathscr{W}(\mathcal{R}(t)) =&~   \displaystyle\int_{\partial \mathcal{R}(t)}\left(\bv^T(\hat{\mathbf{T}}_0-\hat{p}\mathbf{I})-(\hat{\mu} +\alpha \hat{p})\hat{\bh}^v+\dot{\phi} \dfrac{\partial \hat{\Psi}}{\partial \nabla \phi}\right)\cdot \boldsymbol{\nu} ~{\rm d}a,\\
    \mathscr{D}(\mathcal{R}(t)) =&~  \displaystyle\int_{\mathcal{R}(t)}  \left( \hat{\mathbf{T}}_0  + \nabla \phi \otimes \dfrac{\partial \hat{\Psi}}{\partial \nabla \phi} + (\hat{\mu}\phi-\hat{\Psi})\mathbf{I} \right): \nabla \bv  \nn\\
    &\quad \quad \quad -\nabla (\hat{\mu}+\hat{p} \alpha) \cdot \hat{\bh}^v - \dfrac{\partial \hat{\Psi}}{\partial \mathbf{D}} : \dot{\mathbf{D}} -  (\hat{\mu}+\alpha \hat{p})\zeta \hat{\gamma} ~{\rm d}v.\label{eq: def diffusion}
\end{align}
\end{subequations}
Positivity ($\geq$) of the diffusion $\mathscr{D}(\mathcal{R}(t))$ implies, since the control volume $\mathcal{R}(t)$ is arbitrary, positivity of the integrand. To ensure this we impose positivity of the members in the integrand independently. The constitutive class of the free energy $\Psi$, equation \eqref{eq: class Psi}, restricts dependence of $\Psi$ on $\mathbf{D}$. For a more elaborate discussion on the possible dependence of $\Psi$ on $\mathbf{D}$ we refer to  Gurtin \cite{gurtinmodel}. This implies that positivity of the second to last member of the integrand of $\mathscr{D}(\mathcal{R}(t))$ can only be achieved by requiring it to vanish via $\partial \Psi/\partial \mathbf{D} = 0$. This means the constitutive class of $\Psi$ reduces to
\begin{align}\label{eq: class Psi 2}
  \Psi = \hat{\Psi}(\phi,\nabla \phi),
\end{align}
and as a consequence the last member in the integrand in \eqref{eq: def diffusion} vanishes. The energy-dissipation law is thus satisfied when the local inequality is fulfilled:
\begin{align}\label{eq: second law 4}
    \left( \hat{\mathbf{T}}_0  +   \nabla \phi \otimes \dfrac{\partial \hat{ \Psi}}{\partial \nabla \phi} + (\hat{\mu}\phi-\hat{ \Psi})\mathbf{I} \right): \nabla \bv\nn\\
    -\nabla (\hat{\mu}+\alpha \hat{p}) \cdot \hat{\bh}^v-(\hat{\mu}+\alpha \hat{p})\zeta \hat{\gamma} \geq 0.
\end{align}
Based on the inequality \eqref{eq: second law 4}, we now restrict ourselves to stress tensors $\mathbf{T}$, diffusive fluxes $\bh^v$ and mass fluxes $\gamma$ belonging to the constitutive classes:
\begin{subequations}\label{eq: class T h}
\begin{align}
    \mathbf{T} =&~ \hat{\mathbf{T}}(\phi, \nabla \phi, \hat{\mu}, \nabla \hat{\mu}, \mathbf{D}, \hat{p}), \label{eq: class T}\\
    \bh^v =&~ \hat{\bh}^v(\phi, \nabla \phi, \hat{\mu}, \nabla \hat{\mu}, \nabla \hat{p}),\label{eq: class h}\\
    \gamma =&~ \hat{\gamma}(\phi,\hat{\mu},\hat{p}).\label{eq: class gamma}
\end{align}
\end{subequations}
\subsection{Equivalence of alternative modeling restrictions}\label{subsec: Modeling restrictions}
As mentioned in \cref{rmk: arbitrariness of modeling choices}, some of the modeling choices in \cref{sec: dam} seem arbitrary. As a consequence, the modeling restriction obtained in \cref{sec: Derivation of the constitutive modeling restriction} may seem discretionary. In \ref{sec: appendix: alternative constitutive modeling} we provide the derivations of the most obvious alternative modeling restrictions. An overview of the variations is as follows:
\begin{enumerate}
    \item \cref{sec: dam,sec: Derivation of the constitutive modeling restriction}: a volume-measure-based free energy $\Psi$ and the difference of volume fractions $\phi$ as order parameter,
    \item \cref{sec: appendix: alternative constitutive modeling: phi psi}: a mass-measure-based free energy $\psi$ and the difference of volume fractions $\phi$ as order parameter,
    \item \cref{sec: appendix: alternative constitutive modeling: c Psi}: a volume-measure-based free energy $\Psi$ and the difference of concentrations $c$ as order parameter,
    \item \cref{sec: appendix: alternative constitutive modeling: c psi}: a mass-measure-based free energy $\psi$ and the difference of concentrations $c$ as order parameter.
\end{enumerate}
The associated constitutive classes of these modeling choices are:
\begin{subequations}\begin{alignat}{4}
    (1):~&\Psi=\hat{\Psi}(\phi,\nabla \phi), 
    \quad &&\mathbf{T}=\hat{\mathbf{T}}(\phi, \nabla \phi, \hat{\mu}, \nabla \hat{\mu}, \mathbf{D}, \hat{p}), \nn\\
    &\bh^v=\hat{\bh}^v(\phi, \nabla \phi, \hat{\mu}, \nabla \hat{\mu}, \nabla \hat{p}), \quad &&\gamma=\hat{\gamma}(\phi,\hat{\mu},\hat{p}), \\[8pt]
    (2):~&\psi=\doublehat{\psi}(\phi,\nabla \phi), 
    \quad &&\mathbf{T}=\doublehat{\mathbf{T}}(\phi, \nabla \phi, \hat{\mu}, \nabla \hat{\mu}, \mathbf{D}, \doublehat{p}),\nn\\ 
    &\bh^v=\doublehat{\bh}^v(\phi, \nabla \phi, \doublehat{\mu}, \nabla \doublehat{\mu}, \nabla \doublehat{p}),
    \quad &&\gamma=\doublehat{\gamma}(\phi,\doublehat{\mu},\doublehat{p}),\\[8pt]
    (3):~&\Psi=\check{\Psi}(c,\nabla c), 
    \quad &&\mathbf{T}=\check{\mathbf{T}}(c, \nabla c, \check{\mu}, \nabla \check{\mu}, \mathbf{D}, \check{p}),\nn\\
    &\bJ^v=\check{\bJ}^v(c, \nabla c, \check{\mu}, \nabla \check{\mu}, \nabla \check{p}), \quad &&\gamma=\check{\gamma}(c,\check{\mu},\check{p}),\\[8pt]
    (4):~&\psi=\doublecheck{\psi}(c,\nabla c), 
    \quad &&\mathbf{T}=\doublecheck{\mathbf{T}}(c, \nabla c, \doublecheck{\mu}, \nabla \doublecheck{\mu}, \mathbf{D}, \doublecheck{p}),\nn\\
    &\bJ^v=\doublecheck{\bJ}^v(c, \nabla c, \doublecheck{\mu}, \nabla \doublecheck{\mu}, \nabla \doublecheck{p}), \quad &&\gamma=\doublecheck{\gamma}(c,\doublecheck{\mu},\doublecheck{p}),
\end{alignat}
\end{subequations}
where $\hat{p}$, $\doublehat{p}$, $\check{p}$ and $\doublecheck{p}$ are the pressures associated with the modeling choices and the chemical potential-like variables are respectively given by:
\begin{subequations}\label{eq: different chem pot}
  \begin{align}
      (1): \hspace{1cm}\hat{\mu} =&~ \dfrac{ \partial \hat{\Psi}}{\partial \phi} - {\rm div}\dfrac{\partial \hat{\Psi}}{\partial\nabla \phi},\hspace{3cm}\\
      (2): \hspace{1cm}\doublehat{\mu} =&~ \dfrac{ \partial \doublehat{\psi}}{\partial \phi} - \dfrac{1}{\rho}{\rm div}\left(\rho \dfrac{\partial \doublehat{\psi}}{\partial\nabla \phi}\right),\hspace{3cm}\\
      (3): \hspace{1cm}\check{\mu} =&~  \dfrac{ \partial \check{\Psi}}{\partial c} - {\rm div}\dfrac{\partial \check{\Psi}}{\partial\nabla c},\hspace{3cm}\\
      (4): \hspace{1cm}\doublecheck{\mu} =&~ \dfrac{ \partial \doublecheck{\psi}}{\partial c} - \dfrac{1}{\rho}{\rm div}\left(\rho \dfrac{\partial \doublecheck{\psi}}{\partial\nabla c}\right).\hspace{3cm}
  \end{align}
\end{subequations}
Collecting the final restrictions provides:
\begin{subequations}\label{eq: restriction energy-dissipation law}
  \begin{align}
      (1): \hspace{1cm}\left( \hat{\mathbf{T}}_0  +  \nabla \phi\otimes \dfrac{\partial \hat{ \Psi}}{\partial \nabla \phi}  + (\hat{\mu}\phi-\hat{ \Psi})\mathbf{I} \right): \nabla \bv  \hspace{3.5cm}\nonumber \\    \hspace{4cm}
      -\nabla (\hat{\mu}+\alpha \hat{p}) \cdot \hat{\bh}^v -  \zeta (\hat{\mu}+\alpha \hat{p})\hat{\gamma}    \geq 0,\label{eq: restriction energy-dissipation law: phi Psi}\\
      (2): \hspace*{1cm}\left( \doublehat{\mathbf{T}}_0  +  \nabla \phi\otimes \rho\dfrac{\partial \doublehat{ \psi}}{\partial \nabla \phi}  + \rho\doublehat{\mu}\phi\mathbf{I} \right): \nabla \bv\hspace{4.1cm}\nonumber \\    \hspace{4cm}
      -\nabla (\rho\doublehat{\mu} +\alpha \doublehat{p}) \cdot \doublehat{\bh}^v-\zeta(\rho\doublehat{\mu} +\alpha \doublehat{p})\doublehat{\gamma} \geq 0,\label{eq: restriction energy-dissipation law: phi psi}\\
      (3): \hspace{1cm}\left( \check{\mathbf{T}}_0  +  \nabla c \otimes  \dfrac{\partial \check{ \Psi}}{\partial \nabla c} -\check{ \Psi}\mathbf{I} \right): \nabla \bv  \hspace{4.7cm}\nonumber \\    \hspace{4cm}-\nabla \left(\dfrac{\check{\mu}}{\rho}+\beta \check{p}\right) \cdot \check{\bJ}^v - \left(\dfrac{\check{\mu}}{\rho}+\beta \check{p}\right)\check{\gamma} \geq 0, \label{eq: restriction energy-dissipation law: c Psi}\\
      (4): \hspace{1cm}\left( \doublecheck{\mathbf{T}}_0  +\nabla c\otimes  \rho\dfrac{\partial \doublecheck{\psi}}{\partial \nabla c}  \right): \nabla \bv \hspace{5.3cm}\nonumber \\    \hspace{4cm} -\nabla \left(\doublecheck{\mu}+ \doublecheck{p}\beta\right) \cdot \doublecheck{\bJ}^v -\left(\doublecheck{\mu}+ \doublecheck{p}\beta\right)\doublecheck{\gamma}  \geq 0.\label{eq: restriction energy-dissipation law: c psi}
  \end{align}
\end{subequations}
In order to allow comparison of these modeling restrictions, the relation between these modeling classes needs to be specified. We select the following relations:\\

\noindent \textit{free energy classes}:
  \begin{align}\label{eq: identify classes psi}
    \hat{\Psi}(\phi,\nabla\phi) \equiv
    \hat{\rho}(\phi)\doublehat{\psi}(\phi,\nabla\phi)\equiv
    \check{\rho}(c)\doublecheck{\psi}(c,\nabla c)\equiv
    \check{\Psi}(c,\nabla c),
  \end{align}

\noindent \textit{stress tensor classes}:
\begin{align}\label{eq: identify classes stress}
    \hat{\mathbf{T}}(\phi, \nabla \phi, \hat{\mu}, \nabla \hat{\mu},\hat{p}) &~\equiv 
    \doublehat{\mathbf{T}}(\phi, \nabla \phi, \doublehat{\mu}, \nabla \doublehat{\mu},\doublehat{p})  \nn\\
    &~ \equiv 
    \doublecheck{\mathbf{T}}(c, \nabla c, \doublecheck{\mu}, \nabla \doublecheck{\mu},\doublecheck{p})  \nn\\
    &~ \equiv 
    \check{\mathbf{T}}(c, \nabla c, \check{\mu}, \nabla \check{\mu},\check{p}),
\end{align}

\noindent \textit{diffusive flux classes}:

\begin{align}\label{eq: identify classes h j}
    \hat{\bh}^v(\phi, \nabla \phi, \hat{\mu}, \nabla \hat{\mu},\nabla \hat{p}) &~ \equiv
    \hat{\bh}^v(\phi, \nabla \phi, \hat{\mu}, \nabla \hat{\mu},\nabla \doublehat{p}) \nn\\
    &~ \equiv \zeta \check{\bJ}^v(c, \nabla c, \check{\mu}, \nabla \check{\mu},\nabla \check{p})\nn\\
    &~ \equiv\zeta
    \doublecheck{\bJ}^v(c, \nabla c, \doublecheck{\mu}, \nabla \doublecheck{\mu},\nabla \doublecheck{p}),
  \end{align}

\noindent \textit{mass flux classes}:

\begin{align}\label{eq: identify classes mass flux}
    \hat{\gamma}(\phi, \hat{\mu}, \hat{p}) \equiv 
    \doublehat{\gamma}(\phi, \doublehat{\mu}, \doublehat{p})     \equiv 
    \check{\gamma}(c, \check{\mu}, \check{p}) \equiv
    \doublecheck{\gamma}(c, \doublecheck{\mu}, \doublecheck{p}).
\end{align}

The particular identification choice of the diffusive fluxes originates from the relation between the diffusive fluxes before constitutive modeling (\cref{lem: relations h}). Note that a direct consequence of the free energy classes \eqref{eq: identify classes psi} is the equivalence of the associated Korteweg tensors:
  \begin{align}\label{eq: Equivalence of Korteweg tensors}
       \nabla \phi\otimes  \dfrac{\partial \hat{ \Psi}}{\partial \nabla \phi}  = \nabla \phi\otimes \rho\dfrac{\partial \doublehat{ \psi}}{\partial \nabla \phi}  =  \nabla c\otimes  \rho\dfrac{\partial \doublecheck{ \psi}}{\partial \nabla c}  =  \nabla c\otimes \dfrac{\partial \check{ \Psi}}{\partial \nabla c},
  \end{align}
  where we recall the chain rule and the observations $\phi=\phi(c)$ and $c=c(\phi)$.

 \begin{lemma}[Relations between chemical potential-like quantities]\label{lem: Relations chemical potential-like quantities}
The chemical potential-like quantities are related via the following identities:
  \begin{subequations}
    \begin{align}
      \hat{\mu} =&~\rho \doublehat{\mu} +\doublehat{\psi} \jrho,\\
        \check{\mu}=&~ \rho \doublecheck{\mu} -\beta\rho^2\doublecheck{\psi},\\
      \check{\mu} =&~\dfrac{\rho^2}{\rho_1\rho_2} \hat{\mu}.
    \end{align}
  \end{subequations}
\end{lemma}
\begin{proof}
See \ref{appendix: proof of lem and thm}.
\end{proof}

We are now ready to state one of the main results of this paper: \textit{the choice of the order parameter and the type of the free energy do not influence the modeling restrictions}. The selection of a particular order parameter should not be regarded as a modeling step but simply as part of the variable selection as one can easily alter via a variable transformation. 
\begin{theorem}[Equivalence modeling restrictions]\label{thm: mod restrictions}
  Selecting the following relations between the pressures of the various modeling choices:
  \begin{subequations}
    \begin{align}
        \hat{p}=&~\doublehat{p}+\doublehat{\psi} \arho,\\
        \check{p}=&~\doublecheck{p} +\rho\doublecheck{\psi},\\
        \check{p}=&~ \hat{p}+\hat{\mu}\phi,
    \end{align}
  \end{subequations}
  yields equivalence of the restrictions in \eqref{eq: restriction energy-dissipation law}.
\end{theorem}
\begin{proof}
The proof relies on the identifications \eqref{eq: identify classes psi}, \eqref{eq: identify classes stress}, \eqref{eq: identify classes h j} and \eqref{eq: identify classes mass flux}, the relation \eqref{eq: Equivalence of Korteweg tensors} and \cref{lem: Relations chemical potential-like quantities}. Details are provided in \ref{appendix: proof of lem and thm}.
\end{proof}
\subsection{Selection of constitutive models}\label{sec: Selection of constitutive models}
Now that we have established equivalence of the modeling restrictions, we proceed with the one derived in this section:
\begin{align}\label{eq: restate}
    \left( \hat{\mathbf{T}}_0  + \nabla \phi\otimes \dfrac{\partial \hat{ \Psi}}{\partial \nabla \phi}  + (\hat{\mu}\phi-\hat{ \Psi})\mathbf{I} \right): \nabla \bv\hspace{3.5cm}\nonumber \\    \hspace{4cm}-\nabla (\hat{\mu}+\alpha \hat{p}) \cdot \hat{\bh}^v -  \zeta (\hat{\mu}+\alpha \hat{p})\hat{\gamma} \geq 0.
\end{align}
To avoid that variations of $\nabla \bv$ lead to a violation of the modeling restriction, we make the constitutive choice:
\begin{align}\label{eq: stress tensor choice}
    \hat{\mathbf{T}}_0 = - \nabla \phi\otimes \dfrac{\partial \hat{ \Psi}}{\partial \nabla \phi}  - (\hat{\mu}\phi-\hat{ \Psi})\mathbf{I}  + \nu(\phi) (2\mathbf{D}+\lambda({\rm div}\bv) \mathbf{I}),
\end{align}
which is in agreement with \eqref{eq: stress sym} and precludes violation of the energy-dissipation law by the stress contribution. Here $\nu(\phi)\geq 0$ is the dynamic phase-dependent viscosity and $\lambda \geq -2/d$ is a scalar, where we recall that $d$ is the number of dimensions. This form of the viscous stress tensor assumes an isotropic Newtonian mixture. A direct consequence is the expression of the stress tensor $\mathbf{T}$:
\begin{align}\label{eq: stress tensor choice T}
    \hat{\mathbf{T}} = - \nabla \phi\otimes \dfrac{\partial \hat{ \Psi}}{\partial \nabla \phi}  - (\mu\phi-\hat{ \Psi})\mathbf{I}  + \nu(\phi) (2\mathbf{D}+\lambda({\rm div}\bv) \mathbf{I}) -\hat{p}\mathbf{I}.
\end{align}

Next we focus on the diffusive flux $\hat{\bh}^v$. Insisting positivity of the second member in \eqref{eq: restate} implies the form:
\begin{align}\label{eq: model eta h}
    \hat{\bh}^v = -\hat{\mathbf{M}}^v \nabla (\hat{\mu} + \alpha \hat{p}),
\end{align}
for some constitutive quantity $\hat{\mathbf{M}}^v=\hat{\mathbf{M}}^v(\phi,\nabla \phi, \hat{\mu}, \nabla \hat{\mu}, \nabla \hat{p})$, referred to as \textit{mobility tensor}, that is consistent with the inequality
\begin{align}
  -\nabla(\hat{\mu} + \alpha \hat{p})\cdot \left(\hat{\mathbf{M}}^v\nabla(\hat{\mu} + \alpha \hat{p})\right) \leq 0
\end{align}
for all $(\phi,\nabla \phi, \hat{\mu}, \nabla \hat{\mu}, \nabla \hat{p})$.
This may be understood as the \textit{generalized Fick laws of diffusion}.

Lastly, we insist positivity of the third term in \eqref{eq: restate}. Noting that $ \zeta >0$, this implies the form:
\begin{align}\label{eq: const model mass flux}
    \hat{\gamma} = -\hat{m}(\phi,\hat{\mu},\hat{p})  (\hat{\mu}+\alpha \hat{p}),
\end{align}
where the mobility $m=\hat{m}(\phi,\hat{\mu},\hat{p})$ only attains positive values ($m\geq 0)$. Ensuring vanishing mass fluxes in the single fluid region requires $\hat{m}(\phi=\pm1,\hat{\mu},\hat{p})=0$.
\begin{remark}[Constitutive models]
The particular constitutive models of the diffusive flux \eqref{eq: model eta h} and mass flux \eqref{eq: const model mass flux} are the only constitutive models in the specified classes \eqref{eq: class h} and \eqref{eq: class gamma}, that guarantee positivity of the second and third term in \eqref{eq: restate}, respectively. This claim results from a theorem on the solution of thermodynamical inequalities proved by Gurtin in  \cite{gurtin1996generalized}. One could also employ this theorem to construct the most general form of the stress tensor. In this paper we only present a simple constitutive choice of the stress tensor in agreement with unification principle three. 
 \demo
\end{remark}
\begin{remark}[Mobility tensor]
The first appearance of the mobility tensor is in  \cite{gurtin1996generalized}. In phase-field literature it is common to work with an isotropic mobility tensor. \demo
\end{remark}

\begin{remark}[Incompatible mobility in single fluid region]\label{rmk: incomp mob} We show that a degenerate mobility is required in general. Consider in the general non-matching density case ($\alpha \neq 0$) a region $\mathcal{D}_a\in \Omega$ in which $\phi \equiv a$ with $-1\leq a\leq 1$ constant. Note that in $\mathcal{D}_a$ we have $\hat{\Psi}=\hat{\Psi}(\phi,\nabla \phi)=\hat{\Psi}(\phi)=\hat{\Psi}(a)=$~constant, and thus $\hat{\mu}=0$ in $\mathcal{D}_a$. As a consequence we deduce $\hat{\mathbf{M}}^v=\hat{\mathbf{M}}^v(a,\mathbf{0},0,\mathbf{0},\nabla \hat{p})$ in $\mathcal{D}_a$. The mass conservation equation \eqref{eq: mass balance mixture} takes the form:
\begin{align}\label{eq: div v 0}
    {\rm div} \bv = 0,  \quad\quad \text{ in }\mathcal{D}_a,
\end{align}
whereas the phase equation \eqref{eq: phase equations: phi} reduces to:
\begin{align}\label{eq: div v  div h}
   a {\rm div} \bv -\alpha{\rm div}  \left(\hat{\mathbf{M}}^v\nabla  \hat{p}\right) = -\beta \hat{m} \hat{p},  \quad\quad \text{ in }\mathcal{D}_a. 
\end{align}
Therefore:
\begin{align}\label{eq: balance phase field}
    {\rm div}  \left(\hat{\mathbf{M}}^v\nabla  \hat{p}\right) = \zeta \hat{m} \hat{p},  \quad\quad \text{ in }\mathcal{D}_a,
\end{align}
representing a balance equation for the pressure.
Next, consider the cases $a=\pm 1$. In this situation the mass flux contribution vanishes: $m = 0$ and thus \eqref{eq: balance phase field} reduces to:
\begin{align}\label{eq: balance phase field2}
    {\rm div}  \left(\hat{\mathbf{M}}^v\nabla \hat{p}\right) =0,  \quad\quad \text{ in }\mathcal{D}_{\pm1}.
\end{align}
In general \eqref{eq: balance phase field2} represents a balance equation, purely formulated in terms of the pressure, that holds in the single fluid. There exists no such \textit{non-trivial equation} that matches with the standard incompressible Navier-Stokes equations in the single fluid regime. In particular, when $\hat{\mathbf{M}}^v$ is a non-zero constant tensor it follows that:
\begin{align}\label{eq: balance phase field3}
    \Delta \hat{p} =0,  \quad\quad \text{ in }\mathcal{D}_{\pm1},
\end{align}
which in general \textit{does not hold} in the single fluid regime (remark that \eqref{eq: balance phase field3} holds for Stokes flow with a divergence free-force).
Thus in general we are left with the trivial instance of \eqref{eq: balance phase field2}:
\begin{align}\label{eq: mob degenerate}
    \hat{\mathbf{M}}^v = 0,  \quad\quad \text{ in }\mathcal{D}_{\pm1},
\end{align}
i.e. the mobility tensor is of degenerate type.
Note that for many models proposed in literature a non-degenerate mobility was chosen, see \cref{table: overview models constit mod}.\demo
\end{remark}

\begin{remark}[Degenerate mobility]\label{rmk: mob}
The relevance of a degenerate mobility tensor is not new and is in fact well known in community, see e.g. \cite{abels2009diffuse,abels2013incompressible,boyer1999mathematical}. We note that the model presented in Abels et al. \cite{abels2012thermodynamically} has been studied for a variety of mobility choices. Its sharp interface limit has been rigorously shown to exist for both degenerate and non-degenerate mobilities, and the associated sharp interface free boundary problem depends the choice of the mobility. 
\end{remark}
\begin{remark}[Alternative reasoning degenerate mobility]
Instead of following the above arguments leading to a degenerate mobility, one could reason as follows. Note that a substitution of the volume fractions \eqref{eq: order parameters} into the volume-averaged velocity \eqref{eq: volume-averaged velocity} and subsequently into the equation for diffusive flux \eqref{eq: diff flux eq} provides the alternative form of the diffusive flux $\bh^u$:
\begin{align}\label{eq: relation h and jv}
    \bh^u =-\xi(\phi)\jv,
\end{align}
with
\begin{align}
    \xi(\phi) := \phi^2-1 \leq 0.
\end{align}
Obviously $\bh^u$ vanishes whenever $\phi=\pm1$. Recalling \cref{lem: relations h}, we have
\begin{align}\label{eq: relation hv}
    \bh^v&~=\dfrac{\arho}{\rho}\bh^u = -\dfrac{\arho}{\rho}\xi(\phi)\jv,
\end{align} 
and thus the diffusive flux $\bh^v$ also vanishes in the single-fluid case $\phi=\pm1$. Demanding the constitutive model of the diffusive flux $\mathbf{h}^v=\hat{\mathbf{h}}^v$ to be compatible with this restriction we get:
\begin{align}
    \hat{\mathbf{M}}^v\nabla \hat{p} = 0, \quad\quad \text{for } \phi = \pm 1,
\end{align}
where $\hat{\mathbf{M}}^v$ only depends on $\nabla \hat{p})$. As before, the only solution that matches with the standard incompressible Navier-Stokes equations in the single fluid regime is the trivial one:
\begin{align}
    \hat{\mathbf{M}}^v = 0, \quad\quad \text{for } \phi = \pm 1.
\end{align}
Note that \eqref{eq: relation hv} suggests the specific form $\hat{\mathbf{M}}^v = -\xi(\phi)\bar{\mathbf{M}}$ for some $\bar{\mathbf{M}}$.

Of course, this argument can also be made when working with the concentration $c$ as order parameter. On the account of \cref{lem: relations h} and the identity
\begin{align}
    \xi(\phi) = \dfrac{\rho^2}{\rho_1\rho_2}\xi(c),
\end{align}
we have:
\begin{align}
    \bJ^u =-\dfrac{\rho^2\arho}{\rho_1\rho_2}\xi(c)\jv.
\end{align}
Again employing \cref{lem: relations h} reveals that the diffusive flux $\bJ^v$
takes the form:
\begin{align}\bJ^v &~=\dfrac{\rho_1\rho_2}{\arho\rho}\bJ^u = -\rho\xi(c)\jv.
\end{align}
This quantity vanishes in the single-fluid regime $c=\pm1$, and using the constitutive class $\bJ^v=\check{\bJ}^v$ we can deduce that the associated mobility tensor vanishes when $c=\pm1$.
\demo
\end{remark}

\begin{remark}[Alternative diffusive fluxes]
The equivalence of the diffusive fluxes, \eqref{eq: identify classes h j}, provides the form of the diffusive flux $\check{\bJ}^v$:
\begin{align}
    \check{\bJ}^v =&~ -\check{\mathbf{M}}^v\nabla (\hat{\mu}+\alpha \hat{p}),
\end{align}
with $\check{\mathbf{M}}^v=\zeta^{-1}\hat{\mathbf{M}}^v$. Additionally, by identifying the diffusive flux $\hat{\bh}^u = (\rho/\arho)\hat{\bh}^v$ we find:
\begin{align}
    \hat{\bh}^u =&~ -\hat{\mathbf{M}}^u\nabla (\hat{\mu}+\alpha \hat{p}),
\end{align}
with $\hat{\mathbf{M}}^u=(\rho/\arho)\hat{\mathbf{M}}^v$.
\demo
\end{remark}
This completes the constitutive modeling. Via substituting \eqref{eq: stress tensor choice T}, \eqref{eq: model eta h} and \eqref{eq: const model mass flux} into \eqref{eq: model 00} we have obtained \textit{one model} that can be formulated in various ways using:
\begin{itemize}
    \item the mixture velocity $\bv$ or the volume-averaged velocity $\bu$,
    \item the volume fraction difference $\phi$ or the concentration difference $c$,
    \item the mass-measure-based free energy $\hat{\psi}$ or the volume-measure-based free energy $\hat{\Psi}$.
\end{itemize}
We withhold from presenting many formulations and only make variations in the mean velocity to obtain one formulation in terms of the mixture velocity: 
\begin{subequations}\label{eq: model v with const mod}
  \begin{empheq}[left=\empheqlbrace]{align}
   \partial_t (\rho \bv) + {\rm div} \left( \rho \bv\otimes \bv \right) + \nabla \hat{p} + {\rm div} \left( \nabla \phi\otimes \dfrac{\partial \hat{ \Psi}}{\partial \nabla \phi}  + (\mu\phi-\hat{\Psi})\mathbf{I} \right) & \nn\\
    - {\rm div} \left(   \nu (2\mathbf{D}+\lambda({\rm div}\bv) \mathbf{I}) \right)-\rho\mathbf{g} &=~ 0, \label{eq: model v with const mod: mom}\\
 \partial_t \rho + {\rm div}(\rho \bv) &=~ 0, \label{eq: model v with const mod: cont} \\
  \partial_t \phi + {\rm div}(\phi \bv) - {\rm div} \left(\hat{\mathbf{M}}^v\nabla (\hat{\mu}+\alpha \hat{p})\right)  ~+\zeta \hat{m} (\hat{\mu} + \alpha \hat{p})&=~0,\label{eq: model v with const mod: PF}\\
  \hat{\mu} - \dfrac{\partial \hat{\Psi}}{\partial \phi}+{\rm div} \left(  \dfrac{\partial \hat{\Psi}}{\partial \nabla \phi} \right)&=~0,
  \end{empheq}
\end{subequations}
and one formulation in terms of the volume-averaged velocity:
\begin{subequations}\label{eq: model u2 with const mod}
  \begin{empheq}[left=\empheqlbrace]{align}
    \partial_t (\rho \bu + \tilde{\bJ}^u) + {\rm div} \left( \rho \bu\otimes \bu + \tilde{\bJ}^u\otimes \bu + \bu\otimes \tilde{\bJ}^u + \frac{1}{\rho}\tilde{\bJ}^u\otimes\tilde{\bJ}^u \right)&\nn\\
    + \nabla \hat{p} + {\rm div} \left( \nabla \phi\otimes \dfrac{\partial \hat{ \Psi}}{\partial \nabla \phi}  + (\hat{\mu}\phi-\hat{\Psi})\mathbf{I} \right) & \nn\\
    - {\rm div} \left(   \nu (2\nabla^s\left(\bu + \rho^{-1}\tilde{\bJ}^u\right)+\lambda({\rm div}\left(\bu +\rho^{-1}\tilde{\bJ}^u\right)) \mathbf{I}) \right)-\rho\mathbf{g}&=~ 0, \label{eq: model u2 with const mod: mom}\\
 {\rm div} \bu ~+~\beta\hat{m}(\hat{\mu}+\alpha \hat{p})&=~ 0, \label{eq: model u2 with const mod: cont} \\
 \partial_t \phi + \bu\cdot \nabla \phi - {\rm div} \left(\hat{\mathbf{M}}^u\nabla (\hat{\mu}+\alpha \hat{p})\right) +\dfrac{\rho}{2\rho_1\rho_2} \hat{m} (\hat{\mu} + \alpha \hat{p})&~=0\label{eq: model u2 with const mod: PF},\\
 \hat{\mu} - \dfrac{\partial \hat{\Psi}}{\partial \phi}+{\rm div} \left(  \dfrac{\partial \hat{\Psi}}{\partial \nabla \phi} \right)&=~0,
  \end{empheq}
\end{subequations}
with $\tilde{\bJ}^u=-\jrho\hat{\mathbf{M}}^u\nabla(\hat{\mu}+\alpha \hat{p})$. Both models obviously satisfy the exact same form of the energy-dissipation statement \eqref{eq: second law}, which complies with the second unification principle. We remark that the variable transformation that links the models \eqref{eq: model v with const mod} and \eqref{eq: model u2 with const mod} via \eqref{eq: relation mixture velo, volume avg} now involves a constitutive model for the diffusive flux $\tilde{\bJ}^u$.

\section{Unification of existing Navier-Stokes Cahn-Hilliard models}\label{sec: Relation to existing Navier-Stokes Cahn-Hilliard models}

As discussed in the introductory section, \cref{sec: intro}, there exists a wide spectrum with many flavors of NSCH models. It is the purpose of this section to unify them.\\

\noindent \textbf{Lowengrub and Truskinovsky \cite{lowengrub1998quasi}}

To indicate the relation between our model and that of  Lowengrub and Truskinovsky \cite{lowengrub1998quasi} we set $\gamma=0$ and formulate our model in terms of $(\bv, c, \doublecheck{ \psi})$:
\begin{subequations}\label{eq: model c psi v}
  \begin{empheq}[left=\empheqlbrace]{align}
   \partial_t (\rho \bv) + {\rm div} \left( \rho \bv\otimes \bv \right) + \nabla \doublecheck{p} + {\rm div} \left( \nabla c\otimes \rho\dfrac{\partial \doublecheck{ \psi}}{\partial \nabla c}   \right) 
    &\nn\\
    - {\rm div} \left(   \nu(c) (2\mathbf{D}+\lambda({\rm div}\bv) \mathbf{I}) \right)-\rho\mathbf{g} &=~ 0, \label{eq: model c psi v: mom}\\
 \partial_t \rho + {\rm div}(\rho \bv) &=~ 0, \label{eq: model c psi v: cont} \\
  \rho \dot{c} - {\rm div} \left(\hat{\mathbf{M}}^v(\phi(c))\nabla \left(\doublecheck{\mu}+\beta\doublecheck{p} \right) \right) &=~0,\label{eq: model c psi v: PF}\\
  \doublecheck{\mu} - \dfrac{ \partial \doublecheck{ \psi}}{\partial c} + \dfrac{1}{\rho}{\rm div}\left(\rho \dfrac{\partial \doublecheck{ \psi}}{\partial\nabla c}\right)&=~0.\label{eq: model c psi v: mu}
  \end{empheq}
\end{subequations}
To establish \eqref{eq: model c psi v: PF} we have employed the relation between the diffusive classes \eqref{eq: identify classes h j}:
\begin{align}
    \hat{\bh}^v(\phi, \nabla \phi, \hat{\mu}, \nabla \hat{\mu},\nabla \hat{p}) =
    \zeta
    \doublecheck{\bJ}^v(c, \nabla c, \doublecheck{\mu}, \nabla \doublecheck{\mu},\nabla \doublecheck{p})
\end{align}
and the identity
\begin{align}
     \doublecheck{\mu} + \beta \doublecheck{p}=\zeta \left(\hat{\mu} + \alpha \hat{p}\right) .
\end{align}

First we note that the mass balance of the mixture may be written as:
\begin{align}
     {\rm div} \bv + \frac{1}{\rho} \check{\rho}'(c) \dot{c} = 0.
\end{align}
Next, one recognizes the derivative $\partial \rho/\partial c$ in the evolution of the order parameter \eqref{eq: model c psi v: PF} via:
\begin{align}\label{eq: rewritten PF0}
  \rho \dot{c} - {\rm div} \left(\hat{\mathbf{M}}^v(\phi(c))\nabla \left(\doublecheck{\mu}-\dfrac{\doublecheck{p}}{\rho^2} \check{\rho}'(c)\right) \right) =0.
\end{align}
Substitution of \eqref{eq: model c psi v: mu} into \eqref{eq: rewritten PF0} provides the alternative form:
\begin{align}\label{eq: rewritten PF}
  \rho \dot{c} - {\rm div} \left(\hat{\mathbf{M}}^v(\phi(c))\nabla \left(\dfrac{ \partial \doublecheck{ \psi}}{\partial c} - \dfrac{1}{\rho}{\rm div}\left(\rho \dfrac{\partial \doublecheck{ \psi}}{\partial\nabla c}\right)-\dfrac{\doublecheck{p}}{\rho^2} \check{\rho}'(c)\right) \right) =0.
\end{align}
To obtain the model of Lowengrub and Truskinovsky \cite{lowengrub1998quasi} one has to set $\doublecheck{ \psi} = \sigma \varphi(c)/\epsilon + \sigma \epsilon |\nabla c|^2/2$ for some $\varphi=\varphi(c)$, and replace the quantity $\hat{\mathbf{M}}^v(\phi(c))$ in \eqref{eq: rewritten PF} by a non-zero constant isotropic mobility $\tilde{m}\mathbf{I}$:
\begin{subequations}\label{eq: model lowengrub}
  \begin{empheq}[left=\empheqlbrace]{align}
   \partial_t (\rho \bv) + {\rm div} \left( \rho \bv\otimes \bv \right) + \nabla p + \sigma\epsilon {\rm div} \left( \nabla c\otimes \rho \nabla c \right) & \nn\\
    - {\rm div} \left(   \nu (2\mathbf{D}+\lambda({\rm div}\bv) \mathbf{I}) \right) &=~ 0, \label{eq: model lowengrub: mom}\\
 {\rm div} \bv + \frac{1}{\rho} \check{\rho}'(c) \dot{c} &=~ 0, \label{eq: model lowengrub: cont} \\
  \rho \dot{c} - {\rm div} \left(\tilde{m} \nabla \left( \dfrac{\sigma}{\epsilon}\dfrac{\partial \varphi}{\partial c}-\dfrac{\sigma\epsilon}{\rho}{\rm div}\left( \rho \nabla c \right) -\dfrac{p}{\rho^2} \check{\rho}'(c) \right)\right) &=~0.\label{eq: model lowengrub: PF}
  \end{empheq}
\end{subequations}
In other words, this comparison reveals that the quasi-incompressible model of  Lowengrub and Truskinovsky \cite{lowengrub1998quasi} is identical to our model (up to the mobility type, see \cref{rmk: incomp mob}).\\

\noindent \textbf{Shokrpour Roudbari et al. \cite{shokrpour2018diffuse}}, \textbf{Aki et al. \cite{aki2014quasi}} and \textbf{Shen et al. \cite{shen2013mass}}

We now explore the relation of our model with the models \cite{aki2014quasi,shen2013mass,shokrpour2018diffuse}. For the sake of clarity we restrict to $\gamma=0$ and note that only model \cite{aki2014quasi} contains mass fluxes. An alternative formulation of model \eqref{eq: model v with const mod} follows when instead of identity \eqref{eq: id p} one employs
\begin{align}\label{eq: id p2}
    - p^* {\rm div}\bv =&~\frac{p^* \jrho}{\rho} \dot{\phi}\nn\\
                     =&~ \frac{p^*\jrho}{\rho}\left( -{\rm div}\hat{\bh}^v - \phi {\rm div}\bv \right) \nn\\=&~ \nabla \left(\frac{p^*\jrho}{\rho}\right)\cdot\hat{\bh}^v - \frac{p^*\jrho}{\rho} \phi {\rm div}\bv - {\rm div}\left(\frac{p^*\jrho}{\rho}\hat{\bh}^v\right),
\end{align}
in which $p^*$ is the pressure quantity. The associated stress tensor and diffusive flux take the form
\begin{subequations}
\begin{empheq}[left=\empheqlbrace]{align}
    \hat{\mathbf{T}} =&~ - \nabla \phi\otimes \dfrac{\partial \hat{ \Psi}}{\partial \nabla \phi}  - \left(\hat{\mu}\phi-\hat{ \Psi}-\frac{p^*\jrho}{\rho} \phi\right)\mathbf{I} \nn\\
    &~+ \nu(\phi) (2\mathbf{D}+\lambda({\rm div}\bv) \mathbf{I}) -p\mathbf{I},\\
    \hat{\bh}^v =&~ -\hat{\mathbf{M}}^v(\phi)\nabla \left(\hat{\mu}-\frac{p^*\jrho}{\rho} \phi\right).
  \end{empheq}
\end{subequations}
Noting that the mass equation may be written as:
\begin{align}
    {\rm div}  \bv - \alpha {\rm div} \left(\hat{\mathbf{M}}^v(\phi)\nabla \left(\hat{\mu}-\frac{p^*\jrho}{\rho}\right)\right) = 0,
\end{align}
we obtain the equivalent model:
\begin{subequations}\label{eq: model without mu 2}
  \begin{empheq}[left=\empheqlbrace]{align}
   \partial_t (\rho \bv) + {\rm div} \left( \rho \bv\otimes \bv \right) + \nabla p^*& \nn\\
   + {\rm div} \left( \nabla \phi\otimes \dfrac{\partial \hat{ \Psi}}{\partial \nabla \phi}  + \left(\hat{\mu}\phi-\hat{ \Psi}-\frac{p^*\jrho}{\rho} \phi\right)\mathbf{I} \right) & \nn\\
    - {\rm div} \left(   \nu(\phi) (2\mathbf{D}+\lambda({\rm div}\bv) \mathbf{I}) \right)-\rho\mathbf{g} &=~ 0, \label{eq: model S: mom}\\
 {\rm div}  \bv - \alpha {\rm div} \left(\hat{\mathbf{M}}^v(\phi)\nabla \left(\hat{\mu}-\frac{p^*\jrho}{\rho}\right)\right) &=~ 0, \label{eq: model S: cont} \\
  \dot{\phi} +  \phi {\rm div} \bv - {\rm div}  \left(\hat{\mathbf{M}}^v(\phi)\nabla \left(\hat{\mu}-\frac{p^*\jrho}{\rho} \right)\right) &=~0\label{eq: model S: PF}, \\
  \hat{\mu} - \dfrac{\partial \hat{\Psi}}{\partial \phi}+{\rm div}  \left(  \dfrac{\partial \hat{\Psi}}{\partial \nabla \phi} \right)&=~0.
  \end{empheq}
\end{subequations}
The third and fourth terms in the momentum equation \eqref{eq: model S: mom} may be written as:
\begin{align}\label{eq: relation S1}
    &\nabla p^* + {\rm div} \left( \nabla \phi\otimes \dfrac{\partial \hat{ \Psi}}{\partial \nabla \phi}  + \left(\hat{\mu}\phi-\hat{ \Psi}-\frac{p^*\jrho}{\rho} \phi\right)\mathbf{I} \right) = \nn\\
    &\nabla \left(p^*\frac{\arho}{\rho}\right) + {\rm div} \left( \nabla \phi\otimes \dfrac{\partial \hat{ \Psi}}{\partial \nabla \phi}  + \left(\hat{\mu}\phi-\hat{ \Psi}\right)\mathbf{I} \right).
\end{align}
The last member of the phase equation \eqref{eq: model S: PF} may be written as:
\begin{align}\label{eq: relation S2}
    - {\rm div}  \left(\hat{\mathbf{M}}^v(\phi)\nabla \left(\hat{\mu}-\frac{p^*\jrho}{\rho} \right)\right) = - {\rm div}  \left(\hat{\mathbf{M}}^v(\phi)\nabla \left(\hat{\mu}+\alpha p^*\frac{\arho}{\rho} \right)\right).
\end{align}
Equivalence of model \eqref{eq: model without mu 2}
with model \eqref{eq: model v with const mod} follows from \eqref{eq: relation S1} and \eqref{eq: relation S2} via the variable transformation:
\begin{align}
    \hat{p} = p^* \frac{\arho}{\rho} = \frac{p^*}{1-\alpha \phi}.
\end{align} 
Up to the definition of the mobility, this model (and thus also model \eqref{eq: model v with const mod}) is equivalent to the model proposed by  Shokrpour Roudbari et al. \cite{shokrpour2018diffuse}. A variable transformation presented in  Shokrpour Roudbari et al. \cite{shokrpour2018diffuse} reveals that the models \cite{aki2014quasi,shen2013mass} are equivalent, up to the definition of the mobility, to the model  Shokrpour Roudbari et al. \cite{shokrpour2018diffuse}. This reveals that all these models are equivalent to our model (up to the mobility type, see \cref{rmk: incomp mob}).\\

\noindent \textbf{Boyer \cite{boyer2002theoretical}, Ding et al. \cite{ding2007diffuse}} and \textbf{Abels et al. \cite{abels2012thermodynamically}} 

The models \cite{boyer2002theoretical,ding2007diffuse,abels2012thermodynamically} are formulated in terms of the volume-averaged mean velocity. The linear momentum equation deviates from that in our model \eqref{eq: model u2 with const mod}, and as such these models are incompatible with our model. In contrast to the models of Boyer \cite{boyer2002theoretical} and Ding et al. \cite{ding2007diffuse}, the model of Abels et al. \cite{abels2012thermodynamically} is presented with an energy-dissipation law. The kinetic energy in this law is however not an obvious approximation of the kinetic energy of the mixture, see \cref{rmk: Ekin}. Finally, we note that these models are consistent with the incompressible Navier-Stokes equations in the single-fluid regime.

We close this section with an overview of the various models presented in \cref{table: overview models constit mod}.

\begin{table}[h!]
{\small
\begin{tabular}{m{15em}m{5em}m{5em}m{5em}}
\textbf{Model}                      & \rot{\textbf{MT-consistent balance laws}} & \rot{\textbf{Compatible in single fluid}} & \rot{\textbf{Energy-dissipation law}} \\[6pt] \thickhline\\[-4pt]
Abels et al. \cite{abels2012thermodynamically}  & {\color{red}\xmark}  & {\color{darkgreen}\cmark}      & {\color{orange}\smark}              \\[6pt] 
Aki et al. \cite{aki2014quasi}                  & {\color{darkgreen}\cmark}    & {\color{red}\xmark}    & {\color{darkgreen}\cmark}              \\[6pt] 
Boyer \cite{boyer2002theoretical}                    & {\color{red}\xmark}      & {\color{darkgreen}\cmark}   & {\color{red}\xmark}              \\[6pt] 
Ding et al. \cite{ding2007diffuse}                & {\color{red}\xmark}    & {\color{darkgreen}\cmark}   & {\color{red}\xmark}                   \\[6pt] 
Lowengrub and Truskinovsky \cite{lowengrub1998quasi} & {\color{darkgreen}\cmark}     & {\color{red}\xmark} & {\color{darkgreen}\cmark}                  \\[6pt] 
Shen et al. \cite{shen2013mass}               & {\color{darkgreen}\cmark}     & {\color{red}\xmark} & {\color{darkgreen}\cmark}                \\[6pt] 
Shokrpour Roudbari et al. \cite{shokrpour2018diffuse}  & {\color{darkgreen}\cmark}   & {\color{red}\xmark} & {\color{darkgreen}\cmark}              \\[6pt] 
\textit{Current}  & {\color{darkgreen}\cmark}   & {\color{darkgreen}\cmark} & {\color{darkgreen}\cmark}              \\[6pt] \hline
\end{tabular}}
\caption{Comparison of the various NSCH models. The column `MT-consistent balance laws' indicates whether the balance laws of the model are compatible with mixture theory. In the third column `Compatible in single fluid' we state whether the model has a degenerate or non-degenerate mobility, see \cref{rmk: incomp mob}, and in the last column whether the model is energy dissipative. The symbol {\color{orange}\smark} indicates that there is an energy-dissipation law but that the associated kinetic energy is not an obvious approximation of the kinetic energy of the mixture, see also \cref{rmk: Ekin}.}
\label{table: overview models constit mod}
\end{table}

\section{Summary and outlook}\label{sec: discussion}

In this paper we established a unified framework of all existing Navier-Stokes Cahn-Hilliard models. To this purpose, we used the general continuum mixture theory and laid down three unifying principles:
\begin{enumerate}
    \item there is only one system of continuum mechanics balance laws that describes the physical model,
    \item there is only one natural energy-dissipation law that leads to quasi-incompressible Navier-Stokes Cahn-Hilliard models,
    \item variations between the models can only appear in the constitutive choices.
\end{enumerate}
In \cref{sec: mix theory} we provided a precise statement of the principles of mixture theory and their consequences. Furthermore, we showed that the mixture framework leads to one system of balance laws that can be formulated using different variable sets, e.g. in terms of a mass-averaged or volume-averaged velocity. Formulating the balance laws using the volume-averaged velocity, we found a system distinct from existing volume-averaged velocity based models. We illustrated the incompatibility with mixture theory of volume-averaged velocity based models that appear in the literature. This can however easily be repaired. Next, in \cref{sec: 2nd law} we demonstrated how an energy-dissipation law naturally leads to quasi-incompressible Navier-Stokes Cahn-Hilliard models. We proved that the energy-dissipative modeling restriction for the constitutive classes is independent of the variable set. Additionally, we showed that in our framework the mobility tensor of the diffuse flux is of degenerate type. In \cref{sec: Relation to existing Navier-Stokes Cahn-Hilliard models} we demonstrated that, using the appropriate degenerate mobility, existing Navier-Stokes Cahn-Hilliard models are often equivalent reformulations of the same physical model. Finally, we presented an overview of the various Navier-Stokes Cahn-Hilliard models in \cref{table: overview models constit mod}.

While we think that the presented framework and the associated analysis are useful to gain insight into Navier-Stokes Cahn-Hilliard models, we do not claim that they are sufficient in this aspect. We outline some of out thoughts on potential further directions for future research. First, it is essential to establish the sharp interface asymptotics and associated free energy inequalities. We conjecture that the approaches and techniques presented in  Abels et al. \cite{abels2012thermodynamically} and  Aki et al. \cite{aki2014quasi} can directly be applied to the model proposed herein as well. A second important direction for future research lies in the rigorous mathematical analysis of the proposed models. Furthermore, we note that, even though the continuous formulations are equivalent, associated discretization methods are, at least ab initio, not identical. As such, it would be valuable to design numerical schemes that inherit useful properties of the model via convenient formulations. Lastly, it would be worthwhile to explore alternative two-phase flow models that utilize the evolution equation of the diffusive flux.

\appendix

\section{Alternative constitutive modeling}\label{sec: appendix: alternative constitutive modeling}
In this appendix we provide the derivation of the most common alternative constitutive modeling approaches:
\begin{itemize}
    \item \ref{sec: appendix: alternative constitutive modeling: phi psi}: volume fraction $\phi$ and mass-measure-based free energy $\psi$,
    \item \ref{sec: appendix: alternative constitutive modeling: c Psi}: concentration $c$ and volume-measure-based free energy $\Psi$,
    \item \ref{sec: appendix: alternative constitutive modeling: c psi}: concentration $c$ and mass-measure-based free energy $\psi$.
\end{itemize}
\subsection{Constitutive modeling: volume fraction and mass-measure-based free energy}\label{sec: appendix: alternative constitutive modeling: phi psi}
Using the same simplification assumption (S), the total energy $\mathscr{E}$ is the superposition of the Helmholtz free energy, kinetic energy and gravitational energy:
\begin{align}\label{eq: total energy: mass-measure}
  \mathscr{E}(\phi,\bv) =  \displaystyle\int_{\mathcal{R}(t)}(\rho\psi + \mathscr{K}(\mathcal{R}(t)) + \mathscr{G}(\mathcal{R}(t)))~{\rm d}v,
\end{align}
where $\mathcal{R}=\mathcal{R}(t)$ denotes an arbitrary time-dependent control volume in $\Omega$ that is transported by the mixture velocity $\bv$. The mass-measure-based free energy is postulated to pertain to the constitutive class:
\begin{align}\label{eq: class Psi: mass-measure}
  \psi = \doublehat{\psi}(\phi,\nabla \phi,\mathbf{D}),   
\end{align}
and we define the chemical potential-like quantity as:
\begin{align}
    \doublehat{\mu} = \dfrac{ \partial \doublehat{\psi}}{\partial \phi} - \frac{1}{\rho}{\rm div}\left(\rho \dfrac{\partial \doublehat{\psi}}{\partial\nabla \phi}\right).
\end{align}
Analogously to \cref{sec: Derivation of the constitutive modeling restriction}, we work with constitutive classes for stress tensor $\mathbf{T}=\doublehat{\mathbf{T}}$, diffusive flux $\bh^v=\doublehat{\bh}^v$ and mass flux $\gamma=\doublehat{\gamma}$ and postpone their specification.
Next, we postulate the energy-dissipation law:
\begin{align}\label{eq: second law: mass-measure}
    \dfrac{{\rm d}}{{\rm d}t} \mathscr{E}(\phi,\bv) = \mathscr{W}(\mathcal{R}(t)) - \mathscr{D}(\mathcal{R}(t)),
\end{align}
in which $\mathscr{W}(\mathcal{R}(t))$  and $\mathscr{D}(\mathcal{R}(t))\geq 0$ have the same interpretation as in \cref{sec: Derivation of the constitutive modeling restriction}.
With the aid of Reynolds transport theorem and the mixture mass balance law \eqref{eq: mass balance mixture}, the evolution of the free energy in \eqref{eq: total energy: mass-measure} takes the form:
\begin{align}\label{eq: mass measure: evo psi}
      \dfrac{{\rm d}}{{\rm d}t}\displaystyle\int_{\mathcal{R}(t)} \rho \doublehat{\psi} ~{\rm d}v =&~ \displaystyle\int_{\mathcal{R}(t)} \rho \dot{\doublehat{\psi}} ~{\rm d}v\nn\\
      =&~ \displaystyle\int_{\mathcal{R}(t)} \rho \left(\dfrac{\partial \doublehat{\psi}}{\partial \phi} \dot{\phi} + \dfrac{\partial \doublehat{\psi}}{\partial \nabla \phi}\cdot \left(\nabla \phi\right)\dot{}+ \dfrac{\partial \doublehat{\psi}}{\partial \mathbf{D}} : \dot{\mathbf{D}}\right) ~{\rm d}v.
\end{align}
Substitution of \eqref{eq: relation grad phi} into \eqref{eq: mass measure: evo psi} and subsequently the integration by parts provides:
\begin{align}
    \dfrac{{\rm d}}{{\rm d}t}\displaystyle\int_{\mathcal{R}(t)} \rho\doublehat{\psi} ~{\rm d}v  =& \displaystyle\int_{\mathcal{R}(t)} \rho \doublehat{\mu} \dot{\phi} - \left(\nabla \phi\otimes \rho \dfrac{\partial \doublehat{\psi}}{\partial \nabla \phi}\right): \nabla \bv +\rho \dfrac{\partial \doublehat{\psi}}{\partial \mathbf{D}} : \dot{\mathbf{D}}~{\rm d}v  \nn\\
    &~ + \displaystyle\int_{\partial \mathcal{R}(t)}\dot{\phi} \rho\dfrac{\partial \doublehat{\psi}}{\partial \nabla \phi}\cdot \boldsymbol{\nu} ~{\rm d}a.
\end{align}
We now eliminate the  material derivative $\dot{\phi}$ via the substitution of the phase evolution equation \eqref{eq: phase equations: phi}:
\begin{align}\label{eq: psi 3: mass-measure}
    \dfrac{{\rm d}}{{\rm d}t}\displaystyle\int_{\mathcal{R}(t)} \rho\doublehat{\psi} ~{\rm d}v =&~ \displaystyle\int_{\mathcal{R}(t)} \nabla (\rho\doublehat{\mu}) \cdot \doublehat{\bh}^v  - \left(\nabla \phi \otimes\dfrac{\partial \rho\doublehat{\psi}}{\partial \nabla \phi}\right): \nabla \bv + \rho\dfrac{\partial \doublehat{\psi}}{\partial \mathbf{D}} : \dot{\mathbf{D}}  \nn\\
    &\quad\quad\quad-\rho\doublehat{\mu}\phi~{\rm div} \bv+\doublehat{\gamma}\zeta\rho \doublehat{\mu} ~{\rm d}v \nn\\
    &~+ \displaystyle\int_{\partial \mathcal{R}(t)}\left(-\rho\doublehat{\mu} \doublehat{\bh}^v+\dot{\phi}\rho \dfrac{\partial \doublehat{\psi}}{\partial \nabla \phi}\right)\cdot \boldsymbol{\nu} ~{\rm d}a.
\end{align}
Taking now the sum of \eqref{eq: psi 3: mass-measure}, kinetic energy evolution \eqref{eq: kin2} and gravitational energy evolution \eqref{eq: grav2} we arrive at:
\begin{align}\label{eq: second law subst 1 appendix}
    \dfrac{{\rm d}}{{\rm d}t} \mathscr{E}(\phi,\bv) = &~ \displaystyle\int_{\partial \mathcal{R}(t)}\left(\bv^T\doublehat{\mathbf{T}}-\rho\doublehat{\mu} \doublehat{\bh}^v+\dot{\phi} \rho\dfrac{\partial \doublehat{\psi}}{\partial \nabla \phi}\right)\cdot \boldsymbol{\nu} ~{\rm d}a \nn\\
    &~- \displaystyle\int_{\mathcal{R}(t)}  \left( \doublehat{\mathbf{T}}  + \nabla \phi \otimes\rho\dfrac{\partial \doublehat{\psi}}{\partial \nabla \phi} + \mu\phi\mathbf{I} \right): \nabla \bv \nn\\
    &\quad\quad\quad\quad\quad-\nabla (\rho\doublehat{\mu}) \cdot \doublehat{\bh}^v - \dfrac{\partial \doublehat{\psi}}{\partial \mathbf{D}} : \dot{\mathbf{D}} -\zeta(\rho\doublehat{\mu} +\alpha \doublehat{p})\doublehat{\gamma}  ~{\rm d}v.
\end{align}
Employing the partition:
\begin{align}\label{eq: partition psi phi}
    \doublehat{\mathbf{T}} = \doublehat{\mathbf{T}}_0 - \doublehat{p} \mathbf{I} ~\text{with}~\doublehat{\mathbf{T}}_0 := \doublehat{\mathbf{T}} + \doublehat{p} \mathbf{I},
\end{align}
in which $\doublehat{p}$ is a scalar pressure field and the identity \eqref{eq: id p} with pressure $\doublehat{p}$ in \eqref{eq: second law subst 1} allows to identify the rate of work and dissipation as:
\begin{subequations}\label{eq: W, D appendix psi phi}
\begin{align}
    \mathscr{W}(\mathcal{R}(t)) =&~   \displaystyle\int_{\partial \mathcal{R}(t)}\left(\bv^T(\doublehat{\mathbf{T}}_0-\doublehat{p}\mathbf{I})-(\rho\doublehat{\mu} +\alpha \doublehat{p})\doublehat{\bh}^v+\dot{\phi}\rho \dfrac{\partial \doublehat{\psi}}{\partial \nabla \phi}\right)\cdot \boldsymbol{\nu} ~{\rm d}a,\\
    \mathscr{D}(\mathcal{R}(t)) =&~  \displaystyle\int_{\mathcal{R}(t)}  \left( \doublehat{\mathbf{T}}_0  + \nabla \phi \otimes\rho\dfrac{\partial \doublehat{\psi}}{\partial \nabla \phi} + \rho\doublehat{\mu}\phi\mathbf{I} \right): \nabla \bv \nn\\
    &\quad\quad\quad-\nabla (\rho\doublehat{\mu}+\doublehat{p} \alpha) \cdot \doublehat{\bh}^v - \dfrac{\partial \doublehat{\psi}}{\partial \mathbf{D}} : \dot{\mathbf{D}} -\zeta(\rho\doublehat{\mu} +\alpha \doublehat{p})\doublehat{\gamma}  ~{\rm d}v.\label{eq: def diffusion appendix}
\end{align}
\end{subequations}
By following the same argument as in \cref{sec: Derivation of the constitutive modeling restriction} the positivity of the diffusion leads to the reduced class:
\begin{align}\label{eq: class Psi 2: mass-measure}
  \psi = \doublehat{\psi}(\phi,\nabla \phi),
\end{align}
and the modeling restriction readily follows:
\begin{align}\label{eq: second law 3: mass-measure}
    \left( \doublehat{\mathbf{T}}_0  + \nabla \phi\otimes \rho\dfrac{\partial \doublehat{\psi}}{\partial \nabla \phi}  + \rho\doublehat{\mu}\phi\mathbf{I} \right): \nabla \bv \hspace{3cm}\nn\\
    \hspace{3cm}-\nabla (\rho\doublehat{\mu} +\alpha \doublehat{p}) \cdot \doublehat{\bh}^v -\zeta(\rho\doublehat{\mu} +\alpha \doublehat{p})\doublehat{\gamma}\geq 0.
\end{align}
Based on the form of the modeling restriction \eqref{eq: second law 3: mass-measure} we select the following constitutive classes:
\begin{subequations}
\begin{align}
    \mathbf{T} =&~ \doublehat{\mathbf{T}}(\phi, \nabla \phi, \doublehat{\mu}, \nabla \doublehat{\mu}, \mathbf{D}, \doublehat{p}), \label{eq: class T: mass-measure}\\
    \bh^v =&~ \doublehat{\bh}^v(\phi, \nabla \phi, \doublehat{\mu}, \nabla \doublehat{\mu}, \nabla \doublehat{p}),\label{eq: class h: mass-measure}\\
    \gamma =&~ \doublehat{\gamma}(\phi, \doublehat{\mu}, \doublehat{p}).\label{eq: class gamma: mass-measure}
\end{align}
\end{subequations}
\subsection{Constitutive modeling: concentration and volume-measure-based free energy}\label{sec: appendix: alternative constitutive modeling: c Psi}
With the aid of the simplification assumption (S), the total energy $\mathscr{E}$ reads in terms of the concentration $c$:
\begin{align}\label{eq: total energy: Psi c}
  \mathscr{E}(c,\bv) =  \displaystyle\int_{\mathcal{R}(t)}(\Psi + \mathscr{K}(\mathcal{R}(t)) + \mathscr{G}(\mathcal{R}(t)))~{\rm d}v,
\end{align}
where again $\mathcal{R}=\mathcal{R}(t)$ denotes an arbitrary time-dependent control volume in $\Omega$ that is transported by the mixture velocity $\bv$. The volume-measure-based free energy is postulated to belong to the class:
\begin{align}\label{eq: class c Psi}
  \Psi = \check{\Psi}(c,\nabla c,\mathbf{D}),     
\end{align}
and the concentration-based chemical potential-like quantity is the Fr\'{e}chet derivative of the Helmholtz free energy $\Psi$ with respect to $c$:
\begin{align}
    \check{\mu} = \dfrac{ \partial \check{\Psi}}{\partial c} - {\rm div}\dfrac{\partial \check{\Psi}}{\partial\nabla c}.
\end{align}
We work with constitutive classes for stress tensor $\mathbf{T}=\check{\mathbf{T}}$, diffusive flux $\bh^v=\check{\bh}^v$ and mass flux $\gamma=\check{\gamma}$ and postpone their specification.
We postulate the energy-dissipation law:
\begin{align}\label{eq: second law appendix}
    \dfrac{{\rm d}}{{\rm d}t} \mathscr{E}(c,\bv) = \mathscr{W}(\mathcal{R}(t)) - \mathscr{D}(\mathcal{R}(t)),
\end{align}
in which we use the same interpretation of $\mathscr{W}(\mathcal{R}(t))$ and $\mathscr{D}(\mathcal{R}(t))\geq 0$ as in \cref{sec: 2nd law}. To proceed we follow the same procedure as in \cref{sec: 2nd law} with the concentration evolution equation \eqref{eq: phase equations: c} and find:
\begin{align}\label{eq: second law subst 1 c-Psi based}
    \dfrac{{\rm d}}{{\rm d}t} \mathscr{E}(c,\bv) = &~ \displaystyle\int_{\partial \mathcal{R}(t)}\left(\bv^T\check{\mathbf{T}}-\check{\mu} \dfrac{\check{\bJ}^v}{\rho}+\dot{c} \dfrac{\partial \check{\Psi}}{\partial \nabla c}\right)\cdot \boldsymbol{\nu} ~{\rm d}a \nn\\
    &~- \displaystyle\int_{\mathcal{R}(t)}  \left( \check{\mathbf{T}}  + \nabla c\otimes \dfrac{\partial \check{\Psi}}{\partial \nabla c} -\check{\Psi}\mathbf{I} \right): \nabla \bv -\nabla \left(\dfrac{\check{\mu}}{\rho}\right) \cdot \check{\bJ}^v\nn\\
    &\quad\quad\quad\quad- \dfrac{\partial \check{\Psi}}{\partial \mathbf{D}} : \dot{\mathbf{D}} -\dfrac{\check{\mu}}{\rho}\check{\gamma} ~{\rm d}v.
\end{align}
Next, we introduce the partition:
\begin{align}\label{eq: partition appendix Psi c}
    \check{\mathbf{T}} = \check{\mathbf{T}}_0 - \check{p} \mathbf{I} ~\text{with}~\check{\mathbf{T}}_0 := \check{\mathbf{T}} + \check{p} \mathbf{I},
\end{align}
in which $\check{p}$ is a scalar pressure field.  Via the concentration equation \eqref{eq: phase equations: c} we deduce the identity
\begin{align}\label{eq: id p appendix Psi c}
    - \check{p} ~{\rm div}\bv =&~\check{p} \left(  \frac{1}{\rho^2}\frac{\partial \rho}{\partial c}\rho \dot{c} \right) \nn\\
    =&~ -\check{p} \frac{1}{\rho^2}\frac{\partial \rho}{\partial c} \left({\rm div}\check{\bJ}^v-\check{\gamma}\right)\nn\\
     =&~\frac{\check{p}}{2}\left(\frac{1}{\rho_1}-\frac{1}{\rho_2}\right)\left({\rm div}\check{\bJ}^v-\check{\gamma}\right) \nn\\
    =&~\beta \check{p}{\rm div}\check{\bJ}^v-\beta \check{p}\check{\gamma}, 
\end{align}
where we recall the definition:
\begin{align}
    \beta = \dfrac{\rho_2-\rho_1}{2\rho_1\rho_2}= -\dfrac{\jrho}{\rho_1\rho_2}.
\end{align}
Via substitution of the identity \eqref{eq: id p appendix Psi c} we can identify the rate of work and dissipation term:
\begin{subequations}\label{eq: W, D; c based appendix}
\begin{align}
    \mathscr{W}(\mathcal{R}(t)) =&~   \displaystyle\int_{\partial \mathcal{R}(t)}\left(\bv^T\left(\check{\mathbf{T}}_0-\check{p}\mathbf{I}\right)-\left(\dfrac{\check{\mu}}{\rho}+\beta \check{p}\right) \check{\bJ}^v+\dot{c} \dfrac{\partial \check{\Psi}}{\partial \nabla c}\right)\cdot \boldsymbol{\nu} ~{\rm d}a,\\
    \mathscr{D}(\mathcal{R}(t)) =&~  \displaystyle\int_{\mathcal{R}(t)}  \left( \check{\mathbf{T}}_0 + \nabla c\otimes \dfrac{\partial \check{\Psi}}{\partial \nabla c} -\check{\Psi}\mathbf{I} \right): \nabla \bv \nn\\
    &\quad\quad\quad\quad -\nabla \left(\dfrac{\check{\mu}}{\rho}+\beta \check{p}\right) \cdot \check{\bJ}^v - \dfrac{\partial \check{\Psi}}{\partial \mathbf{D}} : \dot{\mathbf{D}} - \left(\dfrac{\check{\mu}}{\rho}+\beta \check{p}\right)\check{\gamma}~{\rm d}v.
\end{align}
\end{subequations}
Using the same argument as before we work with the reduced class of the free energy:
\begin{align}
    \Psi = \check{\Psi}(c,\nabla c),
\end{align}
and obtain the modeling restriction:
\begin{align}\label{eq: second law c}
    \left( \check{\mathbf{T}}_0  + \nabla c\otimes \dfrac{\partial \check{\Psi}}{\partial \nabla c} -\check{\Psi}\mathbf{I} \right): \nabla \bv -\nabla \left(\dfrac{\check{\mu}}{\rho}+\beta \check{p}\right) \cdot \check{\bJ}^v  - \left(\dfrac{\check{\mu}}{\rho}+\beta \check{p}\right)\check{\gamma}\geq 0.
\end{align}
Based on \eqref{eq: second law c} we specify the classes for stress tensor, diffusive flux, and mass flux:
\begin{subequations}\label{eq: appendix class T h: c psi}
\begin{align}
    \mathbf{T} =&~ \check{\mathbf{T}}(c, \nabla c, \mu, \nabla \mu,\mathbf{D}, \check{p}), \label{eq: class T appendix}\\
    \bJ^v =&~ \check{\bJ}^v(c, \nabla c, \mu, \nabla \mu, \nabla \check{p}),\label{eq: class h appendix}    \\
    \gamma =&~ \check{\gamma}(c, \check{\mu}, \check{p}).\label{eq: class gamma appendix}
\end{align}
\end{subequations}
\subsection{Constitutive modeling: concentration and mass-measure-based free energy}\label{sec: appendix: alternative constitutive modeling: c psi}
The simplification assumption (S) provides the 
total energy $\mathscr{E}$ in terms of the 
concentration $c$:
\begin{align}\label{eq: total energy: mass-measure c}
  \mathscr{E}(c,\bv) =  \displaystyle\int_{\mathcal{R}(t)}(\rho\psi + \mathscr{K}(\mathcal{R}(t)) + \mathscr{G}(\mathcal{R}(t)))~{\rm d}v,
\end{align}
where again $\mathcal{R}=\mathcal{R}(t)$ denotes an arbitrary time-dependent control volume in $\Omega$ that is transported by the mixture velocity $\bv$. The mass-measure-based free energy is postulated to belong to the class:
\begin{align}\label{eq: class c psi}
  \psi = \doublecheck{\psi}(c,\nabla c,\mathbf{D}),     
\end{align}
and the concentration-based chemical potential-like quantity is defined as:
\begin{align}
    \doublecheck{\mu} = \dfrac{ \partial \doublecheck{\psi}}{\partial c} - \frac{1}{\rho}{\rm div}\left(\rho \dfrac{\partial \doublecheck{\psi}}{\partial\nabla c}\right).
\end{align}
Analogously to the other cases, we work with constitutive classes for stress tensor $\mathbf{T}=\doublecheck{\mathbf{T}}$, diffusive flux $\bh^v=\doublecheck{\bh}^v$ and mass flux $\gamma=\doublecheck{\gamma}$ and postpone their specification.
We postulate:
\begin{align}\label{eq: second law appendix2}
    \dfrac{{\rm d}}{{\rm d}t} \mathscr{E}(c,\bv) = \mathscr{W}(\mathcal{R}(t)) - \mathscr{D}(\mathcal{R}(t)),
\end{align}
with the same interpretation of $\mathscr{W}(\mathcal{R}(t))$ and $\mathscr{D}(\mathcal{R}(t))\geq 0$. 
We proceed with the evaluation of the evolution of the energy \eqref{eq: total energy: mass-measure c} and find via Reynolds transport theorem and the mass balance equation \eqref{eq: mass balance mixture}:
\begin{align}\label{eq: psi derivation 1: mass-measure c}
    \dfrac{{\rm d}}{{\rm d}t}\displaystyle\int_{\mathcal{R}(t)} \rho\doublecheck{\psi} ~{\rm d}v  =&~ \displaystyle\int_{\mathcal{R}(t)} \rho \dot{\doublecheck{\psi}} ~{\rm d}v\nn\\
    =&~ \displaystyle\int_{\mathcal{R}(t)} \rho \left(\dfrac{\partial \doublecheck{\psi}}{\partial c} \dot{c} + \dfrac{\partial \doublecheck{\psi}}{\partial \nabla c}\cdot \left(\nabla c\right)\dot{}+ \dfrac{\partial \doublecheck{\psi}}{\partial \mathbf{D}} : \dot{\mathbf{D}}\right) ~{\rm d}v.
\end{align}
On the account of the relation
\begin{align}\label{eq: relation grad phi: mass-measure}
    \left(\nabla c\right)\dot{} = \nabla (\dot{c}) - (\nabla c)^T \nabla \bv,
\end{align}
we arrive at:
\begin{align}
    \dfrac{{\rm d}}{{\rm d}t}\displaystyle\int_{\mathcal{R}(t)} \rho\doublecheck{\psi} ~{\rm d}v  =&~ \displaystyle\int_{\mathcal{R}(t)} \rho \doublecheck{\mu} \dot{c} - \left(\nabla c \otimes \dfrac{\partial \rho\doublecheck{\psi}}{\partial \nabla c}\right): \nabla \bv +\rho \dfrac{\partial \doublecheck{\psi}}{\partial \mathbf{D}} : \dot{\mathbf{D}}~{\rm d}v  \nn\\
    &~+ \displaystyle\int_{\partial \mathcal{R}(t)}\dot{c} \rho\dfrac{\partial \doublecheck{\psi}}{\partial \nabla c}\cdot \boldsymbol{\nu} ~{\rm d}a.
\end{align}
By substituting the concentration equation \eqref{eq: phase equations: c} we find:
\begin{align}\label{eq: Psi 2: mass-measure}
    \dfrac{{\rm d}}{{\rm d}t}\displaystyle\int_{\mathcal{R}(t)} \rho\doublecheck{\psi} ~{\rm d}v =&~ \displaystyle\int_{\mathcal{R}(t)} \nabla \doublecheck{\mu} \cdot \doublecheck{\bJ}^v  - \left(\nabla c \otimes \dfrac{\partial \rho\doublecheck{\psi}}{\partial \nabla c}\right): \nabla \bv + \rho\dfrac{\partial \doublecheck{\psi}}{\partial \mathbf{D}} : \dot{\mathbf{D}} + \doublecheck{\mu}\doublecheck{\gamma}   ~{\rm d}v \nn\\
    &~+ \displaystyle\int_{\partial \mathcal{R}(t)}\left(-\doublecheck{\mu} \doublecheck{\bJ}^v+\dot{ c}\rho \dfrac{\partial \doublecheck{\psi}}{\partial \nabla  c}\right)\cdot \boldsymbol{\nu} ~{\rm d}a.
\end{align}
Following the same procedure as in \cref{sec: 2nd law} we find:
\begin{align}\label{eq: second law subst 1 c-psi based}
    \dfrac{{\rm d}}{{\rm d}t} \mathscr{E}(c,\bv) = &~ \displaystyle\int_{\partial \mathcal{R}(t)}\left(\bv^T\doublecheck{\mathbf{T}}-\doublecheck{\mu} \doublecheck{\bJ}^v+\dot{c} \dfrac{\partial \rho\doublecheck{\psi}}{\partial \nabla c}\right)\cdot \boldsymbol{\nu} ~{\rm d}a \nn\\
    &~- \displaystyle\int_{\mathcal{R}(t)}  \left( \doublecheck{\mathbf{T}}  + \nabla c \otimes \dfrac{\partial \rho\doublecheck{\psi}}{\partial \nabla c} \right): \nabla \bv -\nabla \doublecheck{\mu} \cdot \doublecheck{\bJ}^v \nn\\
    &~\quad\quad\quad\quad- \dfrac{\partial \rho\doublecheck{\psi}}{\partial \mathbf{D}} : \dot{\mathbf{D}} -\doublecheck{\mu}\doublecheck{\gamma} ~{\rm d}v.
\end{align}
Via substitution of the partition
\begin{align}\label{eq: partition appendix psi c}
    \doublecheck{\mathbf{T}} = \doublecheck{\mathbf{T}}_0 - \doublecheck{p} \mathbf{I} ~\text{with}~\doublecheck{\mathbf{T}}_0 := \doublecheck{\mathbf{T}} + \doublecheck{p} \mathbf{I},
\end{align}
and the identity \eqref{eq: id p appendix Psi c} with pressure $\doublecheck{p}$ we identify the rate of work and dissipation term:
\begin{subequations}\label{eq: W, D; c based appendix 2}
\begin{align}
    \mathscr{W}(\mathcal{R}(t)) =&~   \displaystyle\int_{\partial \mathcal{R}(t)}\left(\bv^T\left(\doublecheck{\mathbf{T}}_0-\doublecheck{p}\mathbf{I}\right)-\left(\doublecheck{\mu}+\doublecheck{p}\beta\right) \doublecheck{\bJ}^v+\dot{c} \dfrac{\partial \rho \doublecheck{\psi}}{\partial \nabla c}\right)\cdot \boldsymbol{\nu} ~{\rm d}a,\\
    \mathscr{D}(\mathcal{R}(t)) =&~  \displaystyle\int_{\mathcal{R}(t)}  \left( \doublecheck{\mathbf{T}}_0 + \nabla c \otimes  \dfrac{\partial \rho\doublecheck{\psi}}{\partial \nabla c}  \right): \nabla \bv \nn\\
    &\quad\quad\quad\quad -\nabla \left(\doublecheck{\mu}+\doublecheck{p}\beta\right) \cdot \doublecheck{\bJ}^v - \dfrac{\partial \rho\doublecheck{\psi}}{\partial \mathbf{D}} : \dot{\mathbf{D}} -\left(\doublecheck{\mu}+ \doublecheck{p}\beta\right)\doublecheck{\gamma} ~{\rm d}v,
\end{align}
\end{subequations}
where we recall $\beta = -\jrho/(\rho_1\rho_2)$. Using the same argument as before we work with the reduced class of the free energy:
\begin{align}
    \psi = \doublecheck{\psi}(c,\nabla c),
\end{align}
and obtain the modeling restriction
\begin{align}\label{eq: second law c 2}
    \left( \doublecheck{\mathbf{T}}_0  + \nabla c \otimes \rho\dfrac{\partial \doublecheck{\psi}}{\partial \nabla c} \right): \nabla \bv -\nabla \left(\doublecheck{\mu}+\doublecheck{p}\beta\right) \cdot \doublecheck{\bJ}^v -\left(\doublecheck{\mu}+ \doublecheck{p}\beta\right)\doublecheck{\gamma} \geq 0.
\end{align}
Based on the modeling restriction we restrict to the constitutive classes:
\begin{subequations}\label{eq: appendix class T h: c Psi}
\begin{align}
    \mathbf{T} =&~ \doublecheck{\mathbf{T}}(c, \nabla c, \doublecheck{\mu}, \nabla \doublecheck{\mu}, \mathbf{D}, \doublecheck{p}), \label{eq: class T appendix2}\\
    \bJ^v =&~ \doublecheck{\bJ}^v(c, \nabla c, \doublecheck{\mu}, \nabla \doublecheck{\mu}, \nabla \doublecheck{p}),\label{eq: class h appendix2}\\
    \gamma =&~ \doublecheck{\gamma}(c, \doublecheck{\mu}, \doublecheck{p}).\label{eq: class gamma appendix2}
\end{align}
\end{subequations}

\section{Proofs of \cref{lem: Relations chemical potential-like quantities} and \cref{thm: mod restrictions}}\label{appendix: proof of lem and thm}
 \textbf{Lemma 3.1 (Relations between chemical potential-like quantities).} \textit{The chemical potential-like quantities are related via the following identities:}
  \begin{subequations}
    \begin{align}
      \hat{\mu} =&~\rho \doublehat{\mu} +\doublehat{\psi} \jrho,\label{eq: chem 1}\\
        \check{\mu}=&~ \rho \doublecheck{\mu} -\beta\rho^2\doublecheck{\psi},\label{eq: chem 2}\\
      \check{\mu} =&~\dfrac{\rho^2}{\rho_1\rho_2} \hat{\mu}.\label{eq: chem 3}
    \end{align}
  \end{subequations}
\begin{proof}
  The first identity \eqref{eq: chem 1} follows by (i) substituting the first identification from \eqref{eq: identify classes psi} into $\hat{\mu}$, (ii) applying the chain rule and (iii) using the identity $\hat{\rho}'(\phi) = \jrho$.

The second identity \eqref{eq: chem 2} follows in a similar fashion from the third identification in \eqref{eq: identify classes psi} by noting that $ \check{\rho}'(c)=-\beta \rho^2$. 
  
The last identity \eqref{eq: chem 3} is a direct consequence of the chain rule for variational derivatives. Alternatively, one can apply a direct computation which we present here. Substituting the identification \eqref{eq: identify classes psi} into $\check{\mu}$
and expanding the derivatives by the chain rule gives:
  \begin{align}
    \check{\mu} =&~ \dfrac{\partial \hat{\Psi}(\phi(c),\phi'(c)\nabla c)}{\partial c} - \divg \dfrac{\partial \hat{\Psi}(\phi(c),\phi'(c)\nabla c)}{\partial \nabla c} \nn\\
    =&~ \dfrac{\partial \hat{\Psi}(\phi(c),\phi'(c)\nabla c)}{\partial \phi}\phi'(c)+\dfrac{\partial \hat{\Psi}(\phi(c),\phi'(c)\nabla c)}{\partial \nabla\phi}\cdot\nabla c~\phi''(c)\nn\\
    &~- \divg \left(\dfrac{\partial \hat{\Psi}(\phi(c),\phi'(c)\nabla c)}{\partial \nabla \phi}  \right) \phi'(c) -\dfrac{\partial \hat{\Psi}(\phi(c),\phi'(c)\nabla c)}{\partial \nabla\phi}\cdot\nabla (\phi'(c))\nn\\
    =&~ \dfrac{\partial \hat{\Psi}(\phi(c),\phi'(c)\nabla c)}{\partial \phi}\phi'(c)- \divg \left(\dfrac{\partial \hat{\Psi}(\phi(c),\phi'(c)\nabla c)}{\partial \nabla \phi}  \right) \phi'(c) \nn\\
    =&~ \hat{\mu} \phi'(c).
\end{align}
Recalling the relation between $\phi$ and $c$, \eqref{eq: c phi}, completes the proof.
\end{proof}

\textbf{Theorem 3.1 (Modeling restrictions independent of order parameter and free energy type).} \textit{The restrictions in \eqref{eq: restriction energy-dissipation law} are equivalent. In particular we have the relations between the pressures of the various modeling choices:}
  \begin{subequations}
    \begin{align}
        \hat{p}=&~\doublehat{p}+\doublehat{\psi} \arho,\\
        \check{p}=&~\doublecheck{p} +\rho\doublecheck{\psi},\\
        \check{p}=&~ \hat{p}+\hat{\mu}\phi.
    \end{align}
  \end{subequations}
\begin{proof}
  We start off by showing equivalence of the restrictions \eqref{eq: restriction energy-dissipation law: phi Psi} and \eqref{eq: restriction energy-dissipation law: phi psi}. Consider the term in brackets in the first member of \eqref{eq: restriction energy-dissipation law: phi Psi} in isolation. Substituting the identification of the equivalence classes and the identity \eqref{eq: Equivalence of Korteweg tensors}, and applying \cref{lem: Relations chemical potential-like quantities} provides:
  \begin{align}\label{eq: Id stress 1}
      &\hat{\mathbf{T}}_0  + \nabla \phi\otimes \dfrac{\partial \hat{ \Psi}}{\partial \nabla \phi}  + (\hat{\mu}\phi-\hat{ \Psi})\mathbf{I}  =\nn\\
      &\hat{\mathbf{T}} + \nabla \phi\otimes \dfrac{\partial \hat{ \Psi}}{\partial \nabla \phi}  + \left(\hat{p}+\hat{\mu}\phi-\hat{ \Psi}\right)\mathbf{I}  =\nn\\ &\doublehat{\mathbf{T}}  + \nabla \phi\otimes \rho\dfrac{\partial \doublehat{ \psi}}{\partial \nabla \phi}  + \left(\hat{p}+\rho \doublehat{\mu}\phi+\doublehat{\psi} \jrho\phi-\rho\doublehat{ \psi}\right)\mathbf{I}=\nn\\&\doublehat{\mathbf{T}}  + \nabla \phi\otimes \rho\dfrac{\partial \doublehat{ \psi}}{\partial \nabla \phi}  + \rho \doublehat{\mu}\phi\mathbf{I}+\left(\hat{p}-\doublehat{ \psi} \arho\right)\mathbf{I}=\nn\\&\doublehat{\mathbf{T}}_0  + \nabla \phi\otimes \rho\dfrac{\partial \doublehat{ \psi}}{\partial \nabla \phi}  + \rho \doublehat{\mu}\phi\mathbf{I}+\left(\hat{p}-\doublehat{p}-\doublehat{ \psi} \arho\right)\mathbf{I}.
  \end{align}
Next, we consider the term in brackets in the second and third members in isolation and apply \cref{lem: Relations chemical potential-like quantities}:
  \begin{align}\label{eq: Id diff flux 1}
      \hat{\mu}+\alpha \hat{p} &~=\rho \doublehat{\mu}+\doublehat{ \psi} \jrho+\alpha \hat{p}\nn\\
      &~=\rho \doublehat{\mu}+\alpha\left( \hat{p}- \doublehat{ \psi} \arho\right)\nn\\
      &~=\rho \doublehat{\mu}+\alpha \doublehat{p}+\alpha\left( \hat{p}- \doublehat{p}- \doublehat{ \psi} \arho\right).
  \end{align}

Equivalence of \eqref{eq: restriction energy-dissipation law: phi Psi} and \eqref{eq: restriction energy-dissipation law: phi psi} results with the observation that the last members in \eqref{eq: Id stress 1} and \eqref{eq: Id diff flux 1} vanish with the variable transformation:
\begin{align}
    \hat{p}=\doublehat{p}+\doublehat{ \psi} \arho.
\end{align}
We proceed with showing the equivalence of \eqref{eq: restriction energy-dissipation law: c Psi} and \eqref{eq: restriction energy-dissipation law: c psi}. Applying \cref{lem: Relations chemical potential-like quantities} and the identity \eqref{eq: Equivalence of Korteweg tensors} to the term in brackets in the first member of \eqref{eq: restriction energy-dissipation law: c Psi} provides:
\begin{align}\label{eq: Id stress 2}
      \check{\mathbf{T}}_0  + \nabla c\otimes \dfrac{\partial \check{ \Psi}}{\partial \nabla c} -\check{ \Psi}\mathbf{I}  &~=\check{\mathbf{T}}  + \nabla c\otimes \dfrac{\partial \check{ \Psi}}{\partial \nabla c} +\left(\check{p}-\check{ \Psi}\right)\mathbf{I}\nn\\
      &~=\doublecheck{\mathbf{T}}  + \nabla c\otimes \rho\dfrac{\partial \doublecheck{ \psi}}{\partial \nabla c} + \left(\check{p}-\rho\doublecheck{ \psi}\right)\mathbf{I}\nn \\
      &~=\doublecheck{\mathbf{T}}_0  + \nabla \phi\otimes \rho\dfrac{\partial \doublecheck{ \psi}}{\partial \nabla \phi}  +\left(\check{p}-\doublecheck{p}-\rho\doublecheck{ \psi}\right)\mathbf{I}.
  \end{align}
Using \cref{lem: Relations chemical potential-like quantities} the term in brackets in the second and third members of \eqref{eq: restriction energy-dissipation law: c Psi} takes the form:
\begin{align}\label{eq: Id diff flux 2}
    \dfrac{\check{\mu}}{\rho}+\beta\check{p}  &~= \doublecheck{\mu}-\beta\left(\rho\doublecheck{\psi} -\check{p} \right)\nn \\
    &~= \doublecheck{\mu}+\beta\doublecheck{p} +\beta\left(\check{p}-\doublecheck{p} -\rho\doublecheck{\psi} \right). 
\end{align}
Equivalence of \eqref{eq: restriction energy-dissipation law: c Psi} and \eqref{eq: restriction energy-dissipation law: c psi} is a consequence of the observation that the last members in \eqref{eq: Id stress 2} and \eqref{eq: Id diff flux 2} disappear with the variable transformation:
\begin{align}
    \check{p}=\doublecheck{p} +\rho\doublecheck{\psi}.
\end{align}
We finalize the proof by showing equivalence of \eqref{eq: restriction energy-dissipation law: phi Psi} and \eqref{eq: restriction energy-dissipation law: c Psi}. Again on the account of \cref{lem: Relations chemical potential-like quantities} and \eqref{eq: Equivalence of Korteweg tensors}, the term in brackets in the first member of \eqref{eq: restriction energy-dissipation law: phi Psi} can be written as:
\begin{align}\label{eq: Id stress 3}
      &\hat{\mathbf{T}}_0  + \nabla \phi\otimes \dfrac{\partial \hat{ \Psi}}{\partial \nabla \phi}  +\left(\hat{\mu}\phi-\hat{ \Psi}\right)\mathbf{I}  = \nn\\
      &\hat{\mathbf{T}}  + \nabla \phi\otimes \dfrac{\partial \hat{ \Psi}}{\partial \nabla \phi}  +\left(\hat{p}+\hat{\mu}\phi-\hat{ \Psi}\right)\mathbf{I}  = \nn\\&\check{\mathbf{T}}  + \nabla c\otimes \dfrac{\partial \check{ \Psi}}{\partial \nabla c} + \left(\hat{p}+\hat{\mu}\phi-\check{ \Psi}\right)\mathbf{I}=\nn\\&\check{\mathbf{T}}_0  + \nabla c\otimes \dfrac{\partial \check{ \Psi}}{\partial \nabla c} -\check{ \Psi}\mathbf{I}+\left(\hat{p}-\check{p}+\hat{\mu}\phi\right)\mathbf{I}.
  \end{align}
Finally, we focus on the second and third member in \eqref{eq: restriction energy-dissipation law: phi Psi}. On the account of \cref{lem: Relations chemical potential-like quantities} we write the sequence of identities:
\begin{align}\label{eq: Id diff flux 3a}
   & \zeta(\hat{\mu}+\alpha \hat{p}) = \frac{\arho}{\rho_1\rho_2}\hat{\mu} + \beta \hat{p}= \frac{\arho}{\rho^2}\check{\mu}  -\beta \hat{\mu}\phi+\beta \check{p}+\beta(\hat{p}-\check{p}+\hat{\mu}\phi)  \nn\\
   & = \frac{\arho}{\rho^2}\check{\mu}  -\frac{\arho}{\rho^2} \check{\mu}\phi+\beta \check{p}+\beta(\hat{p}-\check{p}+\hat{\mu}\phi)  = \frac{\check{\mu}}{\rho}+\beta \check{p}+\beta(\hat{p}-\check{p}+\hat{\mu}\phi).  
\end{align}
Elimination of the last terms in \eqref{eq: Id stress 3} and \eqref{eq: Id diff flux 3a} via the variable transformation
\begin{align}
    \check{p} = \hat{p}+\hat{\mu}\phi,
\end{align}
and recalling \eqref{eq: identify classes h j} and \eqref{eq: identify classes mass flux} concludes the proof.
\end{proof}

\section*{Acknowledgment}
MtE acknowledges support from the German Research Foundation (Deutsche Forschungsgemeinschaft DFG) via the Walter Benjamin project EI 1210/1-1. KvdZ is grateful to G\"{o}rkem \c{S}im\c{s}ek and Harald van Brummelen for early discussions on the topic. The research by KvdZ was supported by the Engineering and Physical Sciences Research Council (EPSRC), UK, under Grants EP/T005157/1 and EP/W010011/1. DS gratefully acknowledges support from the German Research Foundation (Deutsche Forschungsgemeinschaft DFG) via the Emmy Noether Award SCH 1249/2-1. The authors express their appreciation of the valuable criticism of the anonymous reviewers which has improved the clarity and quality of paper.

\end{document}